\documentclass{article}
        \newcommand{\ul}{\underline}
        \newcommand{\map}{$\longmapsto$ }
        \newcommand{\mm}{$\leftarrow \mapsto$ }
        \newcommand{\cd}{$\Longrightarrow$ }
        \newcommand{\iif}{$\Longleftrightarrow$ }
        \newcommand{\ori }{$\wedge$ }
        \newcommand{\ntt}{$\neg\,$}
        \newcommand{\spc}{\hspace{1cm}}
        \newcommand{\vsp}{\vspace{0.3cm}}
        \newenvironment{proof}{\begin{tabbing}\hspace*{%
        2cm}\=\hspace{9.8cm}\=\hspace{1.5cm}\=\kill}{\end{tabbing}}
\begin{document}
\begin{sffamily}
\LARGE
\begin{center}
REVISED VERSION\\
\normalsize
(Note: A previous version made grandiose unsustainable claims.)\\
\vsp
\vsp
\LARGE
An axiomatic account of space\\ sufficient to provide\\ 
formal proofs of spatial relationships
\end{center} 
\vspace{0.3cm}
\begin{normalsize}
\centerline{ABSTRACT}

\begin{quote}
A definition is given of \emph{seriate sets} as being sets constituted
out of 
structured collections 
of objects which are recursively internally self-similar. Fundamental 
(geometrical) objects of Dimension N are conceived to be constituted out
of 
seriate sets of fundamental objects of Dimension N-1, starting with
points 
assigned to Dimension 0. Syntactical rules to enable such objects to be 
systematically named and combined, are set out. A series of formal
proofs 
of theorems relative to objects of Dimensions 1 and 2 are worked
through. 
A proof that four colours are sufficient to colour any five area map is
given in illustration.
\end{quote}
\end{normalsize}
\vspace{0.3cm}

\Large
INTRODUCTION

\large
This paper describes an axiomatic calculus which can be interpreted as 
descriptive of multidimensional \emph{a-mensurate} space: 
that is space which can be conceived of as existing even before the
notion 
of mensuration has been arrived at. This calculus is founded on the
idea 
of a \emph{seriate set}. Seriate sets are understood to be sets
possessing 
recursive structural internal self-similarity. Seriate sets of points
are 
interpreted as lines, seriate sets of lines as areas, and so on, to 
establish a multidimensional system.

The first section of the paper, entitled \emph{Foundations}, sets out a
purely 
abstract definition of 
seriate sets, that is without any reference to geometrical objects. The 
section following, headed \emph{Dimension 1}, is concerned with seriate 
sets of points as constituting \emph{fundamental lines}. The first two 
theorems here deal with the \emph{external} and \emph{internal} 
self-similairty of fundamental lines: that is with the capacity of two 
such lines to conjoin to form an exactly equivalent fundamental line,
and 
with the ability of any single fundamental line to be divided by any 
(internal) point 
into two exactly equivalent fundamental lines. The remainder of this 
section includes theorems to the effect that any non end-point of a 
fundamental line is between its end-points, while neither of those 
end-points is ever between any other two points of that line; and that
one, 
and one only, of three (different) points within a fundamental line is 
between the other two, for example.

The next section is concerned with objects of Dimension 2. Fundamental 
characteristics of seriate sets 
constituted out of fundamental lines are demonstrated in a number of 
theorems. Such seriate sets are then interpreted as constituting 
fundamental areas. Fundamental areas are shown to be externally and 
internally self-similar analogously to fundamental lines. These results 
then enable a proof to be given that four colours are sufficient for 
any five country map where the countries are all fundamental areas 
within a given fundamental area. 

The paper is divided into separate sections each consisting of a formal 
Definition, or Convention, or Theorem, or the like, 
typically followed by a short section of commentary. A few brief general 
notes are added at the end. An appendix is added giving more rigorous 
proofs of the first two theorems of the main text, which there, for the 
sake of brevity and the reader's convenience, are given in a more
intuitive 
presentation. Illustrative diagrams, purely as an aid to comprehension,
are 
added.

\vsp\Large
FOUNDATIONS

\vsp\large
DEFINITION (D) 0.1: THE AXIOM OF EXISTENCE

The existence of any object can be asserted if that existence does not
entail 
"not A" 
where 'A' stands for any proposition entailed by any characteristic of
any 
array of which the given object is asserted to be a part; and can be 
negatived if it does entail "not A". 

\normalsize
An ARRAY is any collection of 
well-formed objects of the calculus. The assertion of the existence of
an 
object can be imagined to be roughly analogous to drawing a line on a
sheet 
of paper: we can draw any specified line we like so long as this does
not 
offend against something already drawn. For example we can draw a line 
representing the diagonal of a square which does not intersect any other 
line - so long as the other diagonal has not previously been 
drawn. 

\vsp\large
DEFINITION  0.2: SERIATE SETS

A set of more than two like objects constitutes a seriate set S!, if and
only 
if:

(a) no part of any object which is a member of S! is a part of any other
object 
which is also a member of S!;

(b) the membership of S!, includes two, and only two, 
objects, termed "E" objects of S!, while all the other members of S! are 
termed "I" objects of S!, in such wise that given the existence of S!
and 
given any "I" object within S!, the existence can be asserted of a set
of 
seriate sets $\ast$S which fulfils the following conditions:

\begin {itemize}
\item[(i)] the set of all the members of all the seriate sets of $\ast$S
is 
identical one-to-one with all the set of all the members of the given
seriate 
set S!;

\item[(ii)] no member of any given seriate set of $\ast$S is a member of
any 
other seriate set of $\ast$S if it is not an "E" object of that seriate
set;

\item[(iii)] two, and only two, of all the "E" objects of all the
seriate sets 
of $\ast$S, are members of no other seriate set of $\ast$S, and are
identical 
to the two "E" objects of S!;

\item[(iv)] every other "E" object of any seriate set of $\ast$S is a
member 
of two and only two seriate sets of $\ast$S, and the given "I" object is 
identical to such an "E" object;

\item[(v)] no subset of sets of $\ast$S is such that each of the "E" 
objects of every member of that subset is a member of another seriate
set 
of that subset. 
\end{itemize}

\normalsize
See diagrams (at end) for the illustration of a seriate set 
of points in accord with the above.  LIKE objects are all of the same 
dimension, as are points or lines. [Dimension is formally defined p. 8.] 
Variants on the form of S! such as S$^6$! and S$^N$! are all names of 
seriate sets; and variants on the form of $\ast$S such as $\ast$S$^2$
and 
$\ast$S$^J$ are 
all names of (simple) sets, consisting of unstructured aggregations of 
objects. An ELEMENT of an object is any part of that object which 
CONSTITUTES that object, that is, is a MEMBER of the (seriate) set which 
determines the characteristics of the object.

\vsp\large
DEFINITION 0.3: THE AXIOM OF IDENTITY 

One set is IDENTICAL to another if it is not possible to 
distinguish the membership of one set from the membership of the other. 

\normalsize
To DISTINGUISH here refers to the unitary process by which objects are 
isolated/identified/named. 

\vsp\large
DEFINITION 0.4: THE AXIOM OF STABILITY

If, given a seriate set S!, the existence of a seriate set S$^N$! within
S! 
can be asserted, then the existence of any other seriate set within S!
cannot 
be asserted, which set has same "E" objects as S$^N$! and is constituted
by 
a set of objects which is not identical to the set of objects which 
constitutes S$^N$!.

\normalsize
Anything GIVEN it is to be understood to have had its 
existence asserted. An object or set is WITHIN another object or set if 
every element or member of the first object or set is likewise an
element 
or member of the second object or set. D. 0.2 above ensures that a
unique 
serial order can be established between all the objects of a seriate
set. 
It has nothing to say however about any other possible order those same 
objects might (?) have. This axiom ensures that once one serial order in
a 
seriate set has been established no other conflicting order can be
asserted 
to exist.

\vsp\Large
DIMENSION 1

\vsp\large
INTERPRETATION (In.) 1.1

A seriate set of points constitutes a FUNDAMENTAL LINE; and the "E"
objects of 
that seriate set are 
termed the END-POINTS of that line.

\normalsize
A POINT is formally defined in the familiar way as that which has no
parts.

\vsp\large
CONVENTIONS (Cn.) 1.1

(i) The names  of points take the form of unadorned upper case letters,
e.g. 
'A', 'B', 'N'.

(ii) The names of fundamental lines take the form of \ul{AB} 
where 'A' and 'B' are the end-points of the line. The objects - points - 
which constitute that line are referenced by means of a variable, always 
given in lower case, as for instance 'x', in the formulation \ul{AB}(x). 

(iii) The use of round brackets containing a variable immediately 
following the name of an object, or set, is the general syntax for the 
assignment of a variable which references all and only, the elements of
the 
named object or members of the given set. In the case of a fundamental
line 
the name of a point may be included within those same brackets in
addition 
to indicate that it is an element of the given line. E.g. \ul{AB}(x,C) 
states 'C' is a point in the same line. 

(iv) Objects having different 
names are to be assumed non identical (the default condition) until the 
converse is stipulated or proved. 

(v) Symbols and syntax: '\&', '\ori',  stand 
for the logical operators 'and' and 'or'. ' \ntt' is used for negation.
'\cd' 
and '\iif' represent "if....then" and "if and only if"; '\map'
represents a 
one-to-one identity mapping, i.e. 'x \map y' means every object
referenced by 
'x' is identical one-to-one with an object referenced by 'y'; '\mm'
means 
this relationship is reciprocal indicating identity of reference. '+' 
indicates the addition of variables and '-' their subtraction where, and 
only where, the second variable named references a subset of the objects 
referenced by the first. Where the mapping arrow, or the identity sign
'=', 
is used in conjunction with the logical operators '\ori' or '\&': '\ori'
means 
that at least one of the arguments of that operator must have the 
relationship indicated, and '\&' means both must have it. The meaning of
any 
string of symbols in any case where ambiguity might arise is determined
by 
assigning precedence to the leftmost member of any pair of operators.

(vi) In proofs, the writing of the bare unmodified name of an object
represents 
the proposition that that object exists. Such propositions then take
truth 
values. The unadorned name of an object can then appear in any argument 
place of a logical operators, always, it must be remembered, acting as 
surrogate for the proposition that the named object exists.

\normalsize
Sect. (iii): The reference of variables is confined to objects all of
the same 
dimension. 
Sect. (v): [P \iif Q] is more precisely defined as equivalent to 
\{[P \cd Q] \& [Q \cd P]\}. Operator precedence:  \{[A \& B] \map [x
\ori y]\} 
means \{[A \map (x \ori y)] \&  [B \map (x \ori y)]\},  and not \{[(A \&
B) 
\map x]  \ori  [(A \& B) \map y]\}. \{a \map (x \ori y)\} means \{[a
\map x] 
\ori [a \map y]\} and not that some of the objects referenced by 'a' map
onto 
'x' and/or some of them map onto 'y'. \{P \map \ntt a\} is not allowed.
The 
identity sign indicates identity between the objects referenced and not
merely 
identity between signs. Thus  \{C = 
[A \ori B]\} means that some point 'C' is identical to 'A' or 'B' -  or 
possibly both - and not that 'A' or 'B' can be replaced by 'C' in any 
statement containing either.

\vsp\large
THEOREM (Th.) 1.1

Given any two fundamental lines \ul{AB} and \ul{BC} which have no point
in 
common except 'B', then 
the existence of a fundamental line \ul{AC} constituted by all and only 
the points which constitute \ul{AB} and \ul{BC} can be asserted.

\renewcommand{\baselinestretch}{1.6}\small\large
\begin{proof}
1.\>\ul{AB}(x) \& \ul{BC}(y) \\
2.\>\ntt\{[P \map (x \& y)] \& \ntt[P = B]\}\\
3.\>[T \map x] \& \ntt[T = (A \ori B)]\\
4.\>\{\ul{AN}$^1$(n$^1$) \& \ul{N$^1$N}$^2$(n$^2$) \&
\ul{N$^2$N}$^3$(n$^3$) 
\& $\ldots$\\
 \>\>$\ldots$ \& \ul{N$^{K-1}$N}$^K$(n$^K$) \&
\ul{N$^K$B}(n$^{K+1}$)\}\'\\
 \>\>\& \{x \mm [n$^1$ + n$^2$ + n$^3\ldots$ n$^{K+1}$]\}\'\\
 \>\>\& \{T = [N$^1$ \ori N$^2$ \ori N$^3\ldots$ N$^K$]\}\'\>1,3,D.0.2\\
5.\>z \mm (x + y)\\
6.\>[\ul{AN}$^1$(n$^1$) \& \ul{N$^1$N}$^2$(n$^2$) \& $\ldots$ \\
 \>\>$\ldots$ \& \ul{N$^{K-1}$N}$^K$(n$^K$) \& \ul{N$^K$B}(n$^{K+1}$) \& 
 \ul{BC}(y)\}]\'\\
 \>\>\& \{z \mm [n$^1$ + n$^2$ + n$^3\ldots$ n$^{K+1}$ + y]\}\'\\
 \>\>\& \{T = [N$^1$ \ori N$^2$ \ori N$^3\ldots$ N$^K$ \ori
B]\}\'\>1,4,5\\
7.\>\ul{AC}(z)\>\>2,6,It.(T),Sym.D.0.2
\end{proof}

\renewcommand{\baselinestretch}{1}\normalsize
L(LINE) 2:This simply states \ul{AB} and \ul{BC} have no point in common
except 
'B'. 
L. 3. That a point such as 'T',  distinct from 'A' and 'B' can be
asserted 
to map onto 'x' arises from the first clause of D.0.2, but such is 
generally taken to be too obvious to need listing in the third column; 
further it is not generally though necessary to first assert the
existence 
of such an object before asserting such a relationship, except where 
necessary for emphasis. This pattern holds throughout. L. 7: When the 
default condition of Cn.1.1(iv) is read as strongly as possible at L. 4 
then given L. 2, D.0.2(b) clauses (ii)/(iv) hold for 'T' at L. 6. By 
iteration and similarity what follows for 'T' follows for all non 
end-points of \ul{AB} and \ul{BC}. The same holds for 'B' directly from
L. 
1/2. D.0.2(b)v in effect precludes any subset of fundamental lines of
such 
a set as that of L. 4, forming a "separate loop". As it holds for each
of 
\ul{AB} and \ul{BC} in isolation, if it were controverted when they are 
conjoined then they would have to have more than one point in common 
contravening L. 2. All the points of 'z', not being 'A' or 'C' then
fulfil 
the requirements of a seriate set "I" object of \ul{AC}(z). See 
Appendix for a more rigorous proof of this and the next 
theorem.                      

\renewcommand{\baselinestretch}{1}\vsp\large
THEOREM 1.2

Given any non end-point within a fundamental line 
\ul{AB}, then the existence of two fundamental lines \ul{AC} and
\ul{CB}, 
within \ul{AB}, constituted out of all and only the elements of \ul{AB}, 
and having no point in common except C, can be asserted.

\renewcommand{\baselinestretch}{1.6}\small\large
\begin{proof}
1.\>\ul{AB}(C,x)\\
2.\>\{\ul{AX}$^1$(x$^1$) \& \ul{X$^1$X}$^2$(x$^2$) \&
\ul{X$^2$X$^3$}(x$^3$) 
\& $\ldots$\\
 \>\>$\ldots$ \& \ul{X$^{K-1}$X}$^K$(x$^K$) \&
\ul{X$^K$B}(x$^{K+1}$)\}\'\\
 \>\>\& \{x \mm [x$^1$ + x$^2$ + x$^3\ldots$ x$^{K+1}$]\}\'\\
 \>\>\& \{C = [X$^1$ \ori X$^2$ \ori X$^3\ldots$ X$^K$]\}\'\>1,D.0.2\\
3.\>C = X$^1$\\
4.\>\ul{AC}(y$^1$) \& [y$^1$ \mm x$^1$]\>\>2:3\\
5.\>\ul{CX}$^2$(x$^2$) \& \ul{X$^2$X}$^3$(x$^3$) \& $\ldots$ \ul{X$^K$B}
(x$^{K+1}$)\>\>2:3\\
6.\>z$^1$ \mm (x$^2$ + x$^3$ + $\ldots$  x$^{K+1}$)\\
7.\>\ul{CB}(z$^1$) \& [y$^1$ + z$^1$] \mm x\>\>5,6,2,Th.1.1,It\\
8.\>C = X$^2$\\
9.\>\ul{AC}(y$^2$) \& [y$^2$ \mm (x$^1$  + x$^2$)]\>\>8,Sm. 5/7\\
10.\>\ul{CB}(z$^2$) \& [z$^2$ \mm (x$^3$ + x$^4$ + $\ldots$
x$^K$+1)]\>\>
8,Sm 5/7]\\
11.\>[y$^2$ + z$^2$] \mm x\>\>2:9,10\\
12.\>\ul{AC}(y) \& \ul{CB}(z) \& [(y + z) \mm x] \& \\ 
\>\>\ntt\{[P \map (y \& z)] \& \ntt[P = C]\}\'\>2,It.C(X$^N$)
\end{proof}

\renewcommand{\baselinestretch}{1}\normalsize
L. 4: The colon separating the line number in column three indicates the 
embedding of an 
entailment.  L.12: Iterating for 'C' identical to every 'X' in L. 2. D. 
0.2(b) clause (v) obviously holds throughout.  See Appendix for a more 
rigorous proof of this theorem. 

\vsp\large
THEOREM 1.3

Given any two points 'P' and 'Q' within a given fundamental line
\ul{AB}, then 
the existence of a 
fundamental line \ul{PQ} within \ul{AB} together with two other
fundamental 
lines within \ul{AB} - either \ul{AP} and \ul{QB}, or \ul{AQ} and
\ul{PB} - 
can be asserted, such that no point within \ul{AB}, apart from 'P' and
'Q', 
is an element of both \ul{PQ} and either of the those other two
fundamental 
lines, which three lines exhaust \ul{AB}.

\normalsize
This follows directly from Th. 1.2 iterated.

\vsp\large
INTERPRETATION 1.2: BETWEEN

A point 'B' is BETWEEN two other 
points 'A' and 'C', if and only if, within some given set of points, the 
existence of two fundamental lines \ul{AB} and \ul{BC} can be asserted, 
which have in common no other point but 'B'; and this between
relationship 
is said to exist \emph{within} the given set of points. It is symbolised 
A/B/C(x) where 'x' references all the points in the given (underlying) 
set.

\normalsize
It is crucial that throughout one keeps firmly in mind that the 
"between relationship" is a relationship that exists only in respect of 
some underlying set.

\vsp\large
THEOREM 1.4

Every non end-point of a given fundamental line set is between the end-
points 
of the given line, within that line.

\normalsize
This follows directly from Th. 1.2 and In. 1.2.

\vsp\large
THEOREM 1.5

One point of any three (different) points within a fundamental line is
between 
the other two, within that line.

\normalsize
Apply Th. 1.3  to a pair of the given points and then 
Th. 1.2 to the third point, and In. 1.2 will hold in the case of one of 
those given points.

\vsp\large
THEOREM 1.6

If one point is between two others within a given fundamental line, then 
neither of the others is between the other two, within that line.

\renewcommand{\baselinestretch}{1.6}\small\large
\begin{proof}
1.\>\ul{AB}(H,J,K,x)\\
2.\>H/J/K(x)\>\>1,Th.1.5,Rn.\\
3.\>\ul{HJ}(r) \& \ul{JK}(s) \& [(r \& s) \map x]\\
\>\>\& \ntt\{[P \map (r \& s)] \& \ntt[P = J]\}\'\>2,In.1.2\\
4.\>\ul{HK}(t) \& [t \mm (r + s)]\>\>3,Th.1.1\\
5.\>J \map t\>\>3,4\\
6.\>J/H/K(x)\\
7.\>\ul{JH}(r$^1$) \& \ul{HK}(t$^1$) \& [(r$^1$ \& t$^1$) \map x]\\
 \>\>\& \ntt\{[P \map (r$^1$ \& t$^1$)] \& \ntt[P = H]\}\'\>6,In.1.2\\
8.\>\ntt[J \map t$^1$]\>\>7\\
9.\>\ntt[t \mm t$^1$]\>\>5,8\\
10.\>\ntt[H/J/K(x) \& J/H/K(x)]\>\>9,4,7,1,D.0.04
\end{proof}

\renewcommand{\baselinestretch}{1}\normalsize
L. 2: Rename the points H,J,K so J is between the other two if it is not 
already so. It follows from this theorem that 
"between relationships" within a line commutate to form a single
ordering. 
This ordering can then be symbolised through an extrapolation of the
basic 
"between syntax", as for example in A/B/C/D(x).  

\vsp\large
THEOREM 1.7

Given any pair of points within a a given fundamental line then neither 
end-point of that line is between that pair of points, within that line.

\normalsize
Apply Th. 1.2 and this result then follows directly form Th. 1.4 and 
Th. 1.6, in both possible cases.

\vsp\large
INTERPRETATION 1.3 

If, and only if, a set of points constitutes a RING (set), then all and
only, 
those same points, constitute two fundamental lines which have solely in 
common the end-points of each. 

\vsp
THEOREM 1.8

Given any pair of points each within a pair of fundamental lines which 
constitute a ring, then the existence of a pair of 
fundamental lines with those points as end-points can be asserted, which 
pair of lines constitutes a ring constituted out of all and only the
points 
which constitute the original pair of lines.

\renewcommand{\baselinestretch}{1.6}\small\large
\begin{proof}
1.\> \ul{AB}(a) \& \ul{BA}(b)\\
2.\>\ntt\{[P \map (a \& b)] \& \ntt[P = (A \ori B)]\}\\
3.\>[a + b] \mm z\\
4.\>[F \& G] \map a]\\
5.\>\ul{AF}(c) \& \ul{FB}(d) \& [(c + d) \mm a]\>\>4Th.1.2\\
6.\>G \map c\>\>4,5,Rn.\\
7.\>\ul{AG}(e) \& \ul{GF}(f) \& [(e + f) \mm c]\>\>6,Th.1.2\\
8.\>\ul{GF}(g) \& [g \mm (e + b + d)]\>\>7,1,2,5,Th.1.1\\
9.\>[f + g] \mm z\>\>8,7,5,3\\
10.\>\ul{GF}(f) \& \ul{GF}(g) \& [(f + g) \map z]\>\>7/9\\
11.\>\>QED\'10,Sym.,It.
\end{proof}

\renewcommand{\baselinestretch}{1}\normalsize
L. 11: L. 4 assumes the two given points are within the same given line.
By 
symmetry the same result 
holds for both given lines. If the two given points are each within a 
different one of the two given fundamental lines, then by judiciously 
selecting a pair of intermediate points a simple iteration of L.10
achieves 
the same result.

\large
Corollary: Given any pair of points within a ring then the 
existence of a pair of lines with those points as end-points and no
other 
points in common can be asserted constituted out of the same set of
points 
as constitutes the ring.

\normalsize
Given the ring in the first place then the 
situation as described in L.1 must hold by In. 1.3 and hence, given this 
together with the main proof the same applies to every pair of points 
within the ring.

\vsp\large
CONVENTION 1.2

The general form of the name for a ring is X//.../X(x), where 'X' is any
point 
in the ring and 'x' is a variable which 
references all (and only) the points within the set of points that 
constitutes the ring: the names of any points within the ring can be
filled 
in between the slashes, or omitted, as required, in accord with the 
convention for the between relationship.

\vsp
THEOREM 1.9

Every one of three points of a ring is between the other two within the
set 
of the ring.

\renewcommand{\baselinestretch}{1.6}\small\large
\begin{proof}
1.\>A//.../A(x)\\
2.\>[B \& C] \map x\\
3.\>\ul{AB}(a) \& \ul{BA}(b) \& [(a + b) \mm x]\\
 \>\>\& \ntt\{[P  \map (a \& b)] \& \ntt[P = (A \ori
B)]\}\'\>1,2,Th.1.8,Cr.\\
4.\>C \map a\>\>Rn.\\
5.\>\ul{AC}(c) \& \ul{CB}(d) \& [(c + d) \mm a]\\
 \>\>\& \ntt\{[P \map (c \& d)] \& \ntt[P = C]\}\'\>3,4,Th.1.2\\
6.\>\ul{AC}(c) \& \ul{CB}(d) \& \ul{BA}(b) \& [(c \& d \& b) \map
x]\>\>3,5\\
7.\>\ntt[P \map (A \ori B \ori C)] \cd\\               
 \>\>\ntt\{P \map [(c \& d) \ori (c \& b) \ori (d \& b)]\}\'\>3,5\\
8.\>A/C/B(x) \& C/B/A(x) \& B/A/C(x)\>\>6,7,In.1.2.
\end{proof}

\renewcommand{\baselinestretch}{1}\normalsize
L. 4: Reverse the names of the variables 'a' and 'b' if  'C' is not as 
indicated.

\vsp\large
THEOREM 1.10

No set of points which forms a ring is within a fundamental line.

\normalsize
This follows directly from Th. 1.9 and Th. 1.6

\vsp\Large
DIMENSION 2

\vsp\large
INTERPRETATION 2.1

A seriate set constituted out of objects of DIMENSION (N), and the
object 
which that set constitutes, both belong to Dimension (N 
+ 1). A point belongs to Dimension 0. The dimension of a seriate set is 
indicated by a subscript: a set of points, for example, is written 
S$_1$! (the default condition so the subscript is often dropped in this 
case), and a seriate set of lines S$_2$!. When the subscript for 
dimension is used the superscripts differentiating different seriate
sets 
are displaced to the end of the name so that S$_2$!$^N$, and S$_2$!$^M$
are 
names of  different seriate sets of lines.

\normalsize
A set as simply an aggregation 
of objects all of the same the dimension is analogously indicated as
such 
by the use of a subscript indicating the dimension of those objects.
Thus 
$\ast _0$S is a collection of points which has the Dimension 0, as
distinct 
from S$_1$! which indicates a \emph{seriate set} of points which has the 
Dimension 1.                                         

\vsp\large
DEFINITION 2.1 

Given any set of fundamental lines then the set of all (and only) the
points 
within each one of those lines, is termed the 
set of points SUBSUMED by the given set of fundamental lines; and any 
point, or set of such points, within that subsumed set is said to be
WITHIN 
the given set of fundamental lines.                    

\vsp
DEFINITION 2.2

A fundamental line within the set of points subsumed by a given seriate
set of 
fundamental lines S$_2$!, is a SERIATING line with respect to that set,
if and 
only if, no two points within that given fundamental line are both
within the 
same constituting line of S$_2$!.                             

\vsp
CONVENTIONS 2.1

The name of a seriating line takes the form exemplified by \ul{A -
B}(x). 
In any proof (or other context) where only one seriate set of
fundamental lines 
is named, this syntax is taken to indicate that that seriating line has
that 
characteristic with respect to that set.                   

\vsp
DEFINITION 2.3

If, and only if, within a given seriate set, S$_N$! a given object is an
"E" 
object of each 
of two subsidiary seriate sets, S$_N$!$^A$ and S$_N$!$^B$, each within 
S$_N$!, and  S$_N$!$^A$ and S$_N$!$^B$ have no other object in common,
then 
the given object is BETWEEN the other two "E" objects of  S$_N$!$^A$ and 
S$_N$!$^B$, within the given seriate set.                    

\normalsize
This definition simply generalises to all seriate sets the definition of 
between for points in a line given in 
In. 1.2. The syntax remains the same, except that a colon is introduced 
before the brackets containing the variable referencing the objects of
the 
underlying set, where necessary to avoid ambiguity.  E.g. for three
lines 
in a seriate set of lines reference by 'x',  \ul{AB}(a) /  \ul{CD}(b) / 
\ul{EF}(c) : (x)

\vsp\large
CONVENTION 2.2

The name of a seriate set of fundamental lines takes the form 
S$_2$!(x, \ul{S}$_1$!$^1$, \ul{S}$_1$!$^2$), where 'x' references the 
fundamental lines constituting the given seriate set and 
S$_1$!$^1$ and S$_1$!$^2$ reference the "E" objects of that 
set.  

\vsp
THEOREM 2.1 

In a seriate set of fundamental lines no point of any line 
in the set is an element of any other line in the set.

\normalsize
Definition 0.2(a)

\vsp\large
THEOREM 2.2

Seriate sets of fundamental lines are self-similar both 
(a) externally and (b) internally in an exactly analogous way to 
fundamental lines.

Corollary: The between relationships which hold for 
points in a fundamental line hold for the constituting lines of a
seriate 
set of fundamental lines. 

\normalsize
Proofs in Th. 1.1 and 1.2 are unaffected if lines 
replace points in them, and the same follows for the theorems
determining 
between relationships, given Th.  2.1.

\vsp\large
DEFINITION 2.4

A given seriate set of fundamental lines is said to be UNFIXED if it is
not 
part of any array 
within which there exists any object which has within it two points
within 
different constituting line of that given set.

\normalsize
The notion of an unfixed 
seriate set allows us to consider the elements of that set purely in
terms 
of the relationships that arise from their membership of that set, free 
from considerations as to what relationships might be asserted to exist 
between them otherwise.

\vsp\large
THEOREM 2.3 

Given an unfixed seriate set of 
fundamental lines together with any two points on two different 
constituting lines of that set, then the existence of a seriating line
with 
those points as end-points, within the set of points subsumed by the
given 
seriate set, can be asserted.

\renewcommand{\baselinestretch}{1.6}\small\large
\begin{proof}
1.\>S$_2$!(x)\\
2.\>[T \map a] \iif \{[\ul{A$^M$B}$^M$(x$^M$) \map x] \& [T \map
x$^M$]\}\\
3.\>[\ul{A$^P$B}$^P$(x$^P$) \& \ul{A$^Q$B}$^Q$(x$^Q$) \map x])\\
4.\>[P \map x$^P$] \& [Q \map x$^Q$]\\
5.\>S$_2$!$^H$(x$^1$, \ul{S}$_1$!$^P$,\ul{S}$_1$!$^Q$)\\
 \>\>\& [\ul{S}$_1$!$^P$ = \ul{A$^P$B}$^P$(x$^P$)] \& [S$_1$!$^Q$ = 
 \ul{A$^Q$B}$^Q$(x$^Q$)]\'\\
 \>\>\& [x$^1$ \map x]\'\>1,3,Th.2.2\\
6.\>[J \map g]\\
7.\>[J \map x$^N$] \& [\ul{A$^N$B}$^N$(x$^N$) \map x$^1$]\>\>6\\
8.\>\ntt\{[K \map (g \& x$^N$)] \& \ntt[J = K]\}\>\>6,7\\
9.\>[\ul{A$^N$B}$^N$(x$^N$) \map x$^1$] \cd \{J \& [J \map (x$^N$\&
g)]\}\\
10.\>[P \& Q] \map g\\
11.\>$\ast_0$S(g) = \ul{PQ}(g)\>\> 6/10,5,D.0.3\\
12.\>\ul{P-Q}(g) \& [g \map a]\>\>11,9,D.2.2
\end{proof}

\renewcommand{\baselinestretch}{1}\normalsize
L. 2: Here 'a' is defined as a variable 
referencing all the points subsumed by the given seriate set. L. 6/10:
Here 
'g' is defined as a variable referencing only, one and only one point 
within every fundamental line within the seriate set S!$^H$, and
including 
'P' and 'Q'.  L. 11: The set of points reference by 'g' has a one-to-one 
relationship with the seriate set of fundamental lines of S!$^H$,
therefore 
those points cannot be distinguished from those of a seriate set hence
that 
set is identical to a seriate set of points by D.0.3. 

\vsp\large
THEOREM 2.4

Given any unfixed seriate set of fundamental lines then the existence of
a 
ring can 
be asserted constituted out of all and only the end-points of all the 
constituting fundamental lines of the given set together with all the 
points within the "E" objects of the set.

\renewcommand{\baselinestretch}{1.6}\small\large
\begin{proof}
1.\>S$_2$!(x,\ul{S}$_1$!$^1$,\ul{S}$_1$!$^2$)\\
 \>\>\& [S$_1$!$^1$ = \ul{A$^1$B}$^1$(x$^1$)] \& [S$_1$!$^2$ = 
 \ul{A$^2$B}$^2$(x$^2$)]\'\\
2.\>$\ast_0$S(f) \& $\ast_0$S(g) \& \ntt\{P \& [P \map (f \& g)\}\\
3.\>\{[P = (A$^N$ \ori B$^N$)] \& [\ul{A$^N$B}$^N$(x$^N$) \map x]\}\\ 
 \>\>\iif [P \map (f \ori g)]\'\\
4.\>[(A$^1$ \& A$^2$) \map f] \& [(B$^1$ \& B$^2$) \map g]\\
5.\>\ul{A$^1$A}$^2$(f) \& \ul{B$^1$B}$^2$(g)\>\>3/4,Th.2.3,It,D.0.1 \\
6.\>A$^1$/A$^2$/B$^2$/B$^1$/A$^1$(r) \& \{r \map [f + x$^2$ + g +
x$^1$]\}\>\>
1,5,2,Th.1.1,In.1.3
\end{proof}

\renewcommand{\baselinestretch}{1}\normalsize
L. 2/3: Here the two variables 
'f' and 'g' are defined as such that all, and only, the end-points of
all 
the fundamental lines constituting the given seriate set are referenced
by 
one or other of them but not by both. L. 4/5: The end-points of the 
fundamental lines constituting the "E" objects of the given seriate set
are 
here arbitrarily assigned to 'f' and 'g' and then form the end-points of 
two fundamental lines referenced by those variables following the 
\emph{rationale} of Th. 2.3. L. 6: Here one may be inclined to object
that 
a ring cannot necessarily be so constituted because the lines referenced
by 
'f' and 'g' might intersect if the "E" objects of the seriate set were
"the 
wrong way round." No such problem can exist however because the set is 
\emph{unfixed} and hence there is nothing that can define any
relationship 
between the members of the seriate set except membership of that set, so
D. 
0.1 governs: The situation must not be thought of\emph{ }as being one in 
which the lines of the given seriate set  are "already" fixed on a 
page.

\vsp\large
DEFINITION 2.5

The set of points constituting the ring as given 
in Th. 2.4 above is termed the BOUNDARY ((set of) points) of the set of 
points subsumed by the given seriate set.

\vsp
DEFINITION 2.6  

A FIXED seriate set of fundamental lines is identical to a seriate set
of 
fundamental lines 
whose boundary set, with a defined membership, has been asserted to
exist, 
and which is constituted out of the same objects as constitute the fixed 
set.  

\normalsize
The notion of a fixed set is given in this indirect way in order to 
eliminate the need for formal specification of the boundary set and its 
membership, every time the existence of such a set is asserted. An
unfixed 
and a fixed set are not directly converse: for example, if, relative to
a 
given unfixed set, we first assert the existence of a fundamental line 
within that set constituted wholly out of the end-points of the 
constituting lines of the unfixed set, except in that one of its end-
points 
is a non end-point of one of the constituting lines of the given seriate 
set, the unfixed set does not turn  into a fixed set as here defined  
because then the ring so constituted as in Th. 2.4 cannot be asserted to 
exist within the resultant array. 

\vsp\large
THEOREM 2.5

Given a fixed seriate set of fundamental lines together with any two
points on 
two different 
constituting lines of that set, then the existence of a seriating line
with 
those points as end-points, within the set of points subsumed by the
given 
seriate set, can be asserted.

\normalsize
This theorem is the exact repetition of Th. 
2.3 except in so far as the given set is fixed rather than unfixed, and 
holds in the same way as the proof given there, and by D.0.1. This 
identical result which holds for both unfixed and fixed sets as defined 
does not in fact hold In the case of the somewhat unusual array
described 
in the commentary to D. 2.6 above. By considering only unfixed and fixed 
sets, as they have been defined, we eliminate from direct consideration 
that rather anomalous set up.

\vsp\large
THEOREM 2.6

Given any two non intersecting fundamental lines within the boundary set
of a 
given  seriate set of 
fundamental lines, then the existence of another seriate set of
fundamental 
lines can be asserted, which set subsumes the same set of points as is 
subsumed by the given set and has the same boundary set as the first,
and 
whose "E" objects are the two given fundamental lines. 

\renewcommand{\baselinestretch}{1.6}\small\large
\begin{proof}
1.\>S$_2$!$^($x,\ul{S}$_1$!$^1$,\ul{S}$_1$!$^2$)\\
  \>\>\& [\ul{S}$_1$!$^1$ = \ul{A$^1$B}$^1$(x$^1$)] \& [\ul{S}$_1$!$^2$
= 
  \ul{A$^2$B}$^2$(x$^2$)]\'\\
2.\>[P \map a] \iif \{[\ul{A$^N$B}$^N$(x$^N$) \map x] \& [P \map
x$^N$]\}\\
3.\>\{[P = (A$^N$ \ori B$^N$)] \& [\ul{A$^N$B}$^N$(x$^N$) \map x]\}\\ 
 \>\>\iif [P \map (f \ori g)]\'\\
4.\>\ntt\{P \& [P \map (f \& g)]\}\\
5.\>\ul{A$^1$A}$^2$(f) \&  \ul{B$^1$B}$^2$(g)\\
6.\>b \mm [$^f$ + x$^2$ + g  + x$^1$]\\
7.\>\ul{J$^1$K$^1$}(y$^1$) \& \ul{J$^2$K}$^2$(y$^2$) \& [y$^1$ \map
$^f$] \& 
[y$^2$ \map $^g$]\\
 \>\>\& A$^1$/J$^1$/K$^1$(f) \& B$^1$/J$^2$/K$^2$(g)\'\\
8.\>\ul{J$^1$J}$^2$(h) \& [h \map b] \& [x$^1$ \map h]\\
 \>\>\& \ul{K$^1$K}$^2$(w) \& [w \map b] \& [x$^2$ \map w]\'\\
9.\>[T$^1$ \map a] \& \ntt[T$^1$ \map (y$^1$ \ori y$^2$)]\\
10.\>\ul{J$^3$- K}$^3$(y$^3$,T$^1$) \& [y$^3$ \map a] \& [J$^3$ \map h]
\& 
[K$^3$ \map w]\\
 \>\>\& \{\ntt[P = (J$^3$ \ori K$^3$)] \cd \ntt[P \map (y$^3$ \&
b)]\}\'\\
 \>\>\& \ntt\{P \& [P \map y$^3$] \& [P \map (y$^1$ \ori y$^2$)]\}\'\>
 9,Th.2.5,D.0.1\\
11.\>T$^2$ \& [T$^2$ \map a]  \& \ntt[T$^2$ \map (y$^1$ \ori y$^2$ \ori 
y$^3$)]\\
12.\>\ul{J$^4$- K}$^4$(y$^4$,T$^2$) \& [y$^4$ \map a$^]$ \& [J$^4$ \map
h] \& 
[K$^4$ \map w]\\
 \>\>\& \{\ntt[P = (J$^4$ \ori K$^4$)]\} \cd \ntt[P \map (y$^4$ \&
b)]\}\'\\
 \>\>\& \ntt\{P \& [P \map y$^4$] \& [P \map (y$^1$ \ori y$^2$ 
 \ori y$^3$)]\}\'\>11,Th.2.5,D.0.1\\
13.\>$\ast_1$S(y)\\
14.\>\{[(\ul{J$^N$K}$^N$(y$^N$) \& \ul{J$^M$K}$^M$(y$^M$)] \map y] \& 
\ntt[y$^M$ \mm y$^N$]\}\\ 
 \>\>\cd \ntt \{P \& [P \map (y$^N$ \& y$^M$)$^]$\}\'\\
15.\>\{P \& [P \map a]\}\\   
 \>\>\iif \{[\ul{J$^N$K}$^N$(y$^N$) \map y] \& [P \map y$^N$]\}\'\\
16.\>[P \map s] \iif \{[\ul{J$^N$K}$^N$(y$^N$) \map y] \& [P \map
y$^N$]\}\\
17.\>S$_2$!(y, \ul{S}$_1$!$^3$, \ul{S}$_1$!$^4$)]\\
 \>\>\& [S$_1$!$^3$ = \ul{J$^1$K}$^1$(y$^1$)] \& [S$_1$!$^4$ =
\ul{J$^2$K}$^2$
 (y$^2$)]\'\\
 \>\>\& [s \mm a]\'\>13/16,9/12,It.T(a),D.0.3\\
18.\>\>QED\'17,It.
\end{proof}

\renewcommand{\baselinestretch}{1}\normalsize
Note that in the statement of the theorem, because the boundary set is
given 
we do not need 
to state the fixed condition. NON INTERSECTING lines have no point in 
common. Intersecting lines have at least one point in common. L. 1/8 
describe the set-up. L. 10 \& 12: The seriating lines formed as in Th.
2.5 
by taking one point from each of the constituting lines of S!(x) between 
the two constituting lines of that set within which their end-points
lie. 
L. 17: By iterating for every possible 'T' in 'a'  as in L. 9/10,11/12
then 
the existence of a set of non intersecting lines $\ast$(y) can be
asserted 
which set is not distinguishable from a seriate set, hence L. 17. L. 18: 
reference again to L. 9/12 shows the boundary set of the set so formed
to 
be identical to that of the given set. If the given pair of fundamental 
lines are not situated as in L.7 a small number of iterations of L. 17
on 
judiciously selected intermediate pairs of lines achieves the required 
result.

\vsp\large
THEOREM 2.7

Any two seriating line within a fixed seriate set of 
fundamental lines with all their end-points within two of the
constituting 
lines of that set and being non intersecting or, at most having solely
in 
common one end-point, constitute the boundary set of seriate set of 
lines.

\normalsize
This follows from the iteration of Th. 2.2(b) and the use of the 
method of Th. 2.6. Where the seriating lines have an end-point in
common, 
that common point should not be within the fundamental lines selected as 
the "E" objects of the required set.

\vsp\large
THEOREM 2.8

Given any seriate set of fundamental lines then any fundamental line
within the 
set intersects every 
constituting line of that given set which is between any two
constituting 
lines of the set within which lie points of the given line.

\renewcommand{\baselinestretch}{1.6}\small\large
\begin{proof}
1.\>S$_2$!(x,\ul{S}$_1$!$^1$,\ul{S}$_1$!$^2$)\\
2.\>[J \map a] \iif \{[\ul{A$^N$B}$^N$(x$^N$) \map x] \& [J \map
x$^N$]\}\\
3.\>[\ul{S}$_1$!$^1$ = \ul{A$^1$B}$^1$(x$^1$)] \& [\ul{S}$_1$!$^2$ = 
\ul{A$^2$B}$^2$(x$^2$)]\\
4.\>\ul{PQ}(p) \& [p \map a]\\ 
 \>\>\& [P$^1$ \map (p \& x$^P$)] \& [Q$^1$ \map (p \& x$^Q$)]\'\\
 \>\>\& \{[\ul{A$^Q$B}$^Q$(x$^Q$) \& \ul{A$^P$B}$^P$(x$^P$)] \map
x\}\'\\
5.\>\ul{A$^P$B}$^P$(x$^P$) / \ul{A$^N$B}$^N$(x$^N$) /
\ul{A$^Q$B}$^Q$(x$^Q$) 
: (x)\\
6.\>[J \map f] \iif\\ 
 \>\>\{[\ul{A$^1$B}$^1$(x$^1$) / \ul{A$^M$B}$^M$(x$^M$) / 
 \ul{A$^N$B}$^N$(x$^N$): (x)\'\\ 
 \>\>\& (J \map x$^M$)] \ori [J \map x$^1$]\}\'\\
7.\>[J \map g$^]$ \iif\\ 
 \>\>\{[\ul{A$^N$B}$^N$(x$^N$) / \ul{A$^M$B}$^M$(x$^M$) / 
 \ul{A$^2$B}$^2$(x$^2$): (x)\'\\ 
 \>\>\& (J \map x$^M$)] \ori [ J \map x$^2$]\}\'\\
8.\>\ntt\{J \& [J \map (f \& g)]\}\>\>6,7,Th.2.2Cr,1.3\\
9.\>\ntt[T \& [T \map (x$^N$ \& p)]\\
10.\>[K \map p] \cd [K \map (f \ori g)]\>\>9,6,7,1,4\\
11.\>Z \& [Z \map p] \& \ntt[Z \map (f \& g)]\>\>8,10,D.0.2\\
12.\>T \& [T \map (x$^N$ \& p)]\>\>10,11
\end{proof}

\renewcommand{\baselinestretch}{1}\small\large
Corollary: Any fundamental line within a seriate set of fundamental
lines with 
end-points 
one to each of the "E" objects of that set, intersects every
constituting 
line of that set.

\normalsize
L. 1/4 describe the set-up and L. 5 defines a 
constituting line of the seriate set between the two constituting lines
of 
that set within which lie the two given points within the given
fundamental 
line. L. 6,7 define the two sets of all the points "on either side" of
that 
"between" line, which, by Th. 2.2,Cr. and Th. 1.3 applied to every pair
of 
one from each have no point in common (L. 8). L. 9 assumes that the 
"between" line does not intersect the given fundamental line. Thus all
its 
points must lie within the sets of points "to either side" of the
"between" 
line (L. 10). But there must then be a point between every pair of
points 
one from 'f' and one from 'g' within the given line, ergo there must be
a 
point that line not referenced by either, in contradiction of this last
(L. 
11). So there must be a point of the fundamental line within the
"between" 
line (L. 12).

\vsp\large
THEOREM 2.9

Given any fundamental line within a fixed seriate 
set of lines then the existence of a \emph{minimal }set of
\emph{seriately 
catenated} fundamental lines can be asserted, which set of lines
subsumes 
the same set of points as are subsumed by the given lines, and whose 
members are each either wholly within a constituting line of the seriate 
set of lines, or are seriating lines with respect to that set.

\renewcommand{\baselinestretch}{1.6}\small\large
\begin{proof}
1.\>S$_2$!(x,\ul{S}$_1$!$^1$,\ul{S}$_1$!$^2$) \\
2.\>[T \map a] \iif \{[\ul{A$^N$B}$^N$(x$^N$)\map x] \& [T \map
x$^N$]\}\\
3.\>\ul{PQ}(p) \& [p \map a]\\
4.\>[J \map x$^J$] \& [\ul{A$^J$B}$^J$(x$^J$) \map x]\\
5.\>J \map [n$^P$ \ori n$^Q$]\\
6.\>\{[J \map n$^P$] \cd [E = P]\} \& \{[J \map n$^Q$] \cd [E = Q]\}\\
7.\>F \& \ntt\{F$^1$ \& F/F$^1$/J(p) \& \ntt[F$^1$ \map x$^J$]\} 
\& E/F/J(p)\>\>4/6\\
8.\>[J \map n$^B$] \cd [J \map (n$^P$ \& n$^Q$)]\\
9.\>J \map s$^B$\\
10.\> F \& [F \map S$_1$!(x$^F$)] \& G \& [G \map S$_1$!(x$^G$)]\\
 \>\>\& S$_1$!(x$^F$) / S$_1$!(x$^J$) / S$_1$!(x$^G$) : (x)\'\\
 \>\>\& \ntt\{[F$^1$ \& F/F$^1$/J(p) \& (F$^1$ \map x$^J$)]\}\'\\
 \>\>\& \ntt\{[G$^1$ \& G/G$^1$/J(p) \& (G$^1$ \map x$^J$)]\}\'\>9,4\\
11.\>[T \map k] \cd \{[T \map p] \& \ntt[T \map (n$^B$ \ori s$^B$)]\}\\
12.\>[\ul{T$^N$T}$^M$(p$^N$) \map c] \iif\\ 
 \>\>\{[(T$^N$ \& T$^M$) \map k] \& [p$^N$ \map p]\'\\
 \>\>\& \ntt[(T$^{Z}$ \map k) \& T$^M$/T$^{Z}$/T$^N$(p)]\}\'\\
13.\>\ul{Z$^N$Z}$^M$(p$^{Z}$) \map c\\
14.\>\{\ul{Z$^N$Z}$^M$(p$^{Z}$) = \ul{Z$^N$- Z}$^M$(p$^{Z}$)\}\\
 \>\>\ori \{[p$^{Z}$ \map x$^N$] \&  [\ul{A$^N$B}$^N$(x$^N$)\map
x]\}\'\>
 13,12\\
15.\>\>QED\'13/14,12
\end{proof}

\renewcommand{\baselinestretch}{1}\normalsize
A SERIATELY CATENATED set of fundamental lines conforms to $\ast$S in 
D. 0.2(b)ii/iv, and a 
minimal such set is the set with the least number of members which
conforms 
to some other stated condition(s). L. 5/8: n$^B$ points are defined as 
points at which the given line which "continues towards" both 'P' and
'Q' 
within a constituting line of the underlying seriate set. L. 9/10 define 
s$^B$ points as where the given line "traverses" the underlying seriate
set 
"in opposite directions." L. 11: 'k' then references all the points of
the 
given line not referenced by n$^B$ or s$^B$, that is are either its 
end-points or are points where it "varies its traverse" with respect to
the 
underlying set.  L. 12:  'c' here is defined as referencing all the "non 
overlapping" fundamental lines within the given line into which the 'k' 
point "divide it." These lines must then be constituted by "I" points
which 
are either all n$^B$ points and hence are within a constituting line of
the 
given seriate set, or by s$^B$ points and hence are seriating lines with 
respect to it,  because if not, any such line would have "a kink" in it 
contrary to L. 12.

\vsp\large
INTERPRETATION 2.2: FUNDAMENTAL AREAS

A set of points constitutes a FUNDAMENTAL AREA, if and only if, it is
identical 
one-to-one 
to the set of points subsumed by a fixed seriate set of fundamental
lines. 
The name of an area is exemplified by A!(a) where 'a' references all the 
points subsumed by that seriate set of fundamental lines. The addition
of 
an underlined variable within the round brackets, as for example by 'b'
in 
A!(a, \ul{b}), references the boundary (set of) points of the area 
identical to the boundary set relative to which the seriate set of lines
is 
fixed.

\vsp
THEOREM 2.10

Given any two fundamental areas, having in common solely 
a set of boundary points which constitutes a single fundamental line,
then 
the existence of a fundamental area can be asserted which subsumes all
and 
only the points subsumed by the two given fundamental areas, and whose 
boundary set is constituted out of the sum of the boundary sets of the
two 
given areas less the non end-points of the fundamental line common to 
both.

\normalsize
In accord with In. 2.2 and Th. 2.6 the existence can be asserted of 
two seriate sets of fundamental lines, each equivalent to one of the two 
given areas, which seriate sets both have as an "E" object the
fundamental 
line common to the two given areas. By Th. 2.2(a) these sets conjoin to 
form a seriate set of fundamental lines which is equivalent to an area 
subsuming all the points subsumed by the original two areas, and with a 
border as given.

\vsp\large
THEOREM 2.11

Given any fundamental line within a fundamental 
area with solely its end-points within the boundary of the given area,
then 
the existence of two fundamental areas within the given fundamental area 
can be asserted, which areas together subsume the same set of points as
is 
subsumed by the given area and which have no points in common except the 
points within the given fundamental line.

\renewcommand{\baselinestretch}{1.6}\small\large
\begin{proof}
1.\>A!(a$^1$,\ul{b}$^1$)\\
2.\>\ul{JK}(q) \& [q \map a$^1$] \& [(J \& K) \map b$^1$]\\
 \>\>\& \{\ntt[P = (J \ori K)] \cd \ntt[P \map (q \& b$^1$)]\}\'\\
3.\>S$_2$!$^($x,\ul{S}$_1$!$^1$,\ul{S}$_1$!$^2$)\\
 \>\>\& [S$_1$!$^1$ = \ul{A$^1$B}$^1$(x$^1$)] \& [S$_1$!$^2$ = 
 \ul{A$^2$B}$^2$(x$^2$)]\'\\
 \>\>\& [(x$^1$ \& x$^2$) \map b$^1$] \& [J \map x$^1$] \& [K \map
x$^2$]\'\\
4.\>[P \map a$^1$] \iif \{[\ul{A$^N$B}$^N$(x$^N$) \map x] \& 
[P \map x$^N$]\}\>\>1,3,In.2.2,Th.2.6\\
5.\>[P \map f ] \iif \\
 \>\>\{[P = (J$^1$ \ori K$^1$)] \& \ul{J$^1$ - K}$^1$(q$^1$) 
 \& [q$^1$ \map q]  \&\'\\
 \>\>\ntt[\ul{J$^-$ K$^2$}(q$^2$) \& (q$^1$ \map q$^2$) \& 
 \ntt(q$^1$ \mm q$^2$)]\}\'\\
6.\>[Q \map g] \iif\\
 \>\>\{P \& [P \map ($^f$ \& x$^N$)] \& [\ul{A$^N$B}$^N$(x$^N$) \map
x]\'\\
 \>\>\& [Q \map x$^N$] \& [Q \map (q \ori b$^1$)]\}\'\\
7.\>\ul{[J$^1$K}$^1$(q$^1$) \map m] \iif\\
 \>\>\{\ul{[J$^1$- K}$^1$(q$^1$) \& [q$^1$ \map (q \ori b$^1$)] \& 
 [(J$^1$ \& K$^1$) \map g]\'\\
 \>\>\& \ntt[J$^1$/W/K$^1$(q$^1$) \& (W \map g)]\}\'\\
8.\>[\ul{N$^1$N}$^2$(n$^1$) \map n] \cd\\
 \>\>\{[\ul{A$^N$B}$^N$(x$^N$) \map x] \& [n$^1$ \map x$^N$]\}\'\\
9.\>A!(a$^S$, \ul{r}$^S$) \map s\\
10.\>[\ul{M$^1$N}$^1$(m$^1$) \& \ul{M$^2$N}$^2$(m$^2$)] \map m\\
11.\>[\ul{M$^1$M}$^2$(n$^1$) \map n] \& \{[\ul{N$^1$N}$^2$(n$^2$) 
\map n] \ori [N$^1$ = N$^2$]\}\\
12.\>\ntt\{[\ul{M$^3$N}$^3$(m$^3$) \map m] \& [M$^3$ \map n$^1$] \& 
[N$^3$ \map n$^2$]\\                
 \>\>\& M$^1$/M$^3$/M$^2$(n$^1$)\}\'\\
13.\>\{\ntt[N$^1$ = N$^2$] \cd [r$^S$ \mm [m$^1$ + n$^2$ + m$^2$ +
n$^1$]\\
 \>\>\& \{[N$^1$ = N$^2$] \cd [r$^S$ \mm [m$^1$ + m$^2$ + n$^1$]\}\'\\
 \>\>\& \{M$^1$/N$^1$/N$^2$/M$^2$/M$^1$(r$^S$) \ori  
 M$^1$/N$^1$/M$^2$/M$^1$(r$^S$)\}\'\\
 \>\>\& [a$^S$ \map a$^1$]\'\>9,Th.1.1,Th.2.7\\
14.\>[\ul{M$^1$N}$^1$(m$^1$) \map m] \& [m$^1$ \map q]\\
15.\>\{[A!(a$^N$, \ul{r}$^N$) \& A!(a$^M$, \ul{r}$^M$)] \map s\}\\
 \>\>\& [m$^1$ \map (r$^N$ \& r$^M$)] \& \ntt[r$^N$ \mm r$^M$]\'\\
 \>\>\& [A!(a$^N$, \ul{r}$^N$) \map c] \& [A!(a$^M$, \ul{r}$^M$) 
 \map d]\'\>14 \\
16.\>\{[A!(a$^N$, \ul{r}$^N$) \& A!(a$^M$, \ul{r}$^M$)] \map s\} 
\& \ntt[a$^M$ \mm a$^N$]             \\
 \>\>\& [\ul{M$^1$M}$^2$(n$^M$) \map n] \& [n$^M$ \map (r$^N$ \&
r$^M$)]\'\\
17.\>\{[A!(a$^N$, \ul{r}$^N$) \map c] \cd [A!(a$^M$, \ul{r}$^M$) \map
c]\}\\
 \>\>\& \{[A!(a$^N$, \ul{r}$^N$) \map d] \cd [A!(a$^M$, \ul{r}$^M$) \map
d]\}
 \'\>16\\
18.\>[P \map a$^2$] \iif \{[P \map a$^N$] \& [A!(a$^N$, \ul{r}$^N$) \map
c]\}\\
19.\>[P \map a$^3$] \iif \{[P \map a$^N$] \& [A!(a$^N$, \ul{r}$^N$) \map
d]\}\\
20.\>A!(a$^2$,\ul{b}$^2$) \& A!(a$^3$,\ul{b}$^3$)\\
 \>\>\& \{[P \map (b$^2$ \ori b$^3$)] \cd [P \map (q \ori b$^1$)]\}\'\\
 \>\>\& \ntt\{P \& \ntt[P \map q] \& [P \map (a$^2$ \& a$^3$)]\}\'\\
 \>\>\& [a$^1$ \mm (a$^2$ + a$^3$)]\'\>Th.2.9,Th. 2.10,It.
\end{proof}

\renewcommand{\baselinestretch}{1}\normalsize
This theorem together with the last show how fundamental areas are
externally 
and internally self-similar. 
L. 5: 'f' references the end-points of all the seriating lines of a
minimal 
seriate catenation  equivalent to the given line. L. 6: 'g' 
references\emph{every} intersection point of every constituting line of 
the given seriate set having an 'f' point within it, with the given
line, 
or the with boundary of the given area. L. 7 defines all the seriating 
lines within the given line or the boundary which have 'g' points as 
end-points and no similar point between them. L. 9/13 define as 's'
areas 
all the areas within the given area defined by any two 'm' lines with 
end-points on the same two constituting lines of the underlying seriate 
set, and with no such line between them. L. 14/15: Each of the two 's' 
areas "lying on opposite sides" of every seriating section of the given 
line are then assigned each to one of two mutually exclusive classes,
those 
areas "lying side-by-side" being assigned to the same class (L. 16/17).
All 
the areas in each one of these classes are then added together in accord 
with Th. 2.10 to form two areas with Th. 2.9 ensuring they are together 
equivalent to the original area. If the consistency of the assignment of 
areas to 'c' and 'd' is doubted, application of Th. 2.8 shows the
agreement 
between L. 15 \& L. 17.

\vsp\large
THEOREM 2.12

Given a ring within a fundamental area then the existence can be
asserted of 
one and only one fundamental area 
within the given fundamental area with that ring as its boundary (set). 

\normalsize
The existence of the area with the given ring as boundary follows from a 
maximum of two, judicious partitionings of the given fundamental area in 
accord with Th. 2.11. If another area with the same boundary set could 
exist within the same given area then the assignment of points by Th.
2.11 
would be controverted.

\vsp\large
THEOREM 2.13

Given any pair of fundamental areas having in common only boundary
points, 
together with a fundamental line 
wholly within those areas, and having a non boundary point of each
within 
it, then a point within that fundamental line must also be a common 
boundary point of those areas.

\normalsize
If the two areas have merely one point in 
common then the given fundamental line must have that point within it or
it 
can not be such by D. 0.2. Otherwise, iterated judicious distinction of 
seriate sets equivalent to the given areas and application of Th. 2.8 
governs.

\vsp\large
THEOREM 2.14

If three fundamental lines within a given fundamental 
area have solely in common a pair of points which are the end-points of
all 
of them, then one and only one pair of those three fundamental lines
forms 
the boundary set of a fundamental area whose existence within the given 
fundamental area can be asserted, which area subsumes all the points
within 
the other two fundamental areas whose existence can be asserted within
the 
given area, defined by boundary sets constituted out of the other two
pairs 
of lines from the original three. 

\renewcommand{\baselinestretch}{1.6}\small\large
\begin{proof}
1.\>A!(a)\\
2.\>\ul{PQ}(x) \& \ul{PQ}(y) \& \ul{PQ}(z) \& [(x \& y \& z) \map a]\\
3.\>\ntt[T = (P \ori Q)] \\
4.\>\ntt[T \map [(x \& y) \ori (x \& z) \ori (y \& z)]\>\>3\\
5.\>A!(a$^1$,\ul{b}$^1$) \& [a$^1$ \map a$^]$ \& [b$^1$ \mm (x + y)]\>\>
1,2,Th.2.12\\
6.\>A$^!$(a$^2$,\ul{b}$^2$) \& [a$^2$ \map a] \& [b$^2$ \mm (x + z)]\>\>
1,2,Th.2.12\\
7.\>A$^!$(a$^3$,\ul{b}$^3$) \& [a$^3$ \map a] \& [b$^3$ \mm (z + y)]\>\>
1,2,Th.2.12\\
8.\>[X \map x] \& [Y \map y] \& [Z \map z]\\  
 \>\>\&  \ntt[(X \ori Y \ori Z) = (P \ori Q)]\'\\
9.\>[Z \map a$^1$]\\
10.\>[(a$^2$ \& a$^3$) \map a$^1$] \& [Z \map (a$^1$ \& a$^2$ \&
a$^3$)]\>\>
9, Th.2.13 \\
11.\>\ntt[X \map a$^3$] \& \ntt[Y \map a$^2$]\>\>9,Th.2.12,2.11\\
12.\>\ntt[Z \map a$^1$] \cd \{[X \map a$^3$] \ori [Y \map a$^2$]\}\>\>
Th.2.12,2.11\\
13.\>\>QED\'9/11,11,Sym.
\end{proof}

\renewcommand{\baselinestretch}{1}\normalsize
L. 10: As 'Z' a non end-point of the 'z' 
line is within a$^1$ then every point of that line must be within a$^1$
or 
else Th. 2.13 is controverted by L. 3/4 and similarly every point of
a$^2$ 
must be within a$^1$. L. 11: Th. 2.11 and Th. 2.12 require that given L. 
9/10 a unique layout arises with 'X' and 'Y' in different subsidiary
areas. 
L. 12: If 'Z' is not in a$^1$ then, if 'X' is not in a$^3$, then add
a$^1$ 
to a$^3$ to form a$^2$ and 'Y' must be within it, and this covers all 
possibilities.

\vsp\large
THEOREM 2.15

No three fundamental areas all within a given fundamental area, which
have 
solely boundary points in common, can all have 
the same one point in common, which point is a non end-point of a 
fundamental line constituted out of boundary points of two of them.

\normalsize
Assume such a point exists. Then consider a ring within the conjunction
of the 
two 
areas which have the fundamental line in common, with the postulated
point 
as a non boundary point of the area within that ring. If a third area,
not 
within the other two, has that point within it then by Th. 2.5 the 
existence of a fundamental line can be asserted having within it the 
postulated point and also a point not within that ring, which line does
not 
intersect the ring, in contravention of Th. 2.13. 

\vsp\Large
A Five Country Map

\large
Five fundamental areas, each within a given fundamental area, such 
that they have no points in common except boundary points, cannot have
each 
with every other a common set of points which constitutes a fundamental 
line.

\renewcommand{\baselinestretch}{1.6}\small\large
\begin{proof}
1.\>A!(a)\\
2.\>\ul{P$^N$P}$^M$(j$^)$ \map t\\
3.\>A!(a$^N$,\ul{b}$^N$) \& A!(a$^M$,\ul{b}$^M$)\\
 \>\spc \& \ntt[a$^M$ \mm a$^N$] \& [(a$^N$ \& a$^M$) \map a]\\
 \>\spc \& \{\ntt[Z \map (b$^N$ \& b$^M$)] \cd \ntt[Z \map (a$^N$ \&
a$^M$)]\}
 \\
 \>\spc \& \ul{FG}(k) \& [k \map (b$^N$ \& b$^M$)]\\
4.\>[P$^N$ \map a$^N$] \& \ntt[P$^N$ \map b$^N$]\\
  \>\spc \& [P$^M$ \map a$^M$] \& \ntt[P$^M$ \map b$^M$]\\
  \>\spc \& [j \map (a$^N$ + a$^M$)]\\
  \>\spc \& [T \map (j \& k)] \& \ntt[T = (F  \ori G)]\\
 \>\spc \& \{\ntt[T$^1$ = T] \cd \ntt[T$^1$ \map (j \& k)]\}\\
5.\>\{[\ul{P$^N$P}$^1$(j$^1$) \& \ul{P$^N$P}$^2$(j$^2$)] \map t\} \&
\ntt[j$^1$ 
\mm j$^2$]\\
6.\>\ntt[T = P$^N$] \cd \ntt[T \map (j$^1$ \& j$^2$)]\>\>5,D.0.1\\
7.\>\{ \ul{P$^A$P}$^C$(j$^A$C) \& \ul{P$^C$P}$^B$(j$^C$B)]\\
 \>\spc \& \ul{P$^A$P}$^D$(j$^A$D)\\
 \>\spc \& \ul{P$^D$P}$^B$(j$^D$B)\\
 \>\spc \& \ul{P$^A$P}$^E$(j$^A$E)\\
 \>\spc \& \ul{P$^E$P}$^B$(j$^E$B)\\
 \>\spc \& \ul{P$^A$P}$^B$(j$^A$B)\} \map t\\
8.\>\ul{P$^A$P}$^B$(w$^C$) \& \ul{P$^A$P}$^B$(w$^D$) \&
\ul{P$^A$P}$^B$(w$^E$)
\\
 \>\spc \& [w$^C$ \mm (j$^A$C + j$^C$B)]\\
 \>\spc \& [w$^D$ \mm (j$^A$D + j$^D$B)]\\
 \>\spc \& [w$^E$ \mm (j$^A$E + j$^E$B)]\\
 \>\spc \&  [w \mm j$^A$B]\>\>7,Th.1.1\\
9.\>$\ldots$ /P$^A$//P$^B$/$\ldots$(h$^1$) \& [h$^1$ \mm (w$^C$ +
w$^E$)]\\
 \>\spc \& A!(f$^1$,\ul{h}$^1$) \& [f$^1$ \map a]\\
 \>\spc \& [P$^D$ \map a$^1$]\>\> 5/6,8,Th.2.14,Rn.\\
10.\>$\ldots$ /P$^A$//P$^B$/$\ldots$(h$^2$) \& [h$^2$ \mm (w + w$^C$)]\\
 \>\spc \& A!(f$^2$,\ul{h}$^2$) \& [f$^2$ \map a]\\
 \>\spc \& [P$^E$ \map f$^2$]\>\>5/6,8,Th.2.14,Rn.\\
11.\>\ntt\{\ul{P$^C$P}$^E$ \& [\ul{P$^C$P}$^E$ \map
t]\}\>\>10,Th.2.13,2.15\\
12.\>$\ldots$/P$^A$//P$^B$/$\ldots$(h$^3$) \& [h$^3$ \mm (w + w$^E$)]\\
 \>\spc \& A!(f$^3$,\ul{h}$^3$) \& [f$^3$ \map a$^1$]\\
 \>\spc \& [P$^D$ \map f$^3$]\>\>5/6,8,Th.2.14,Rn.\\
13.\>\ntt\{\ul{P$^C$P}$^D$ \& [\ul{P$^C$P}$^D$ \map
t]\}\>\>12,Th.2.13,2.1\\
14.\>\>QED\'11,13,Sym.
\end{proof}

\renewcommand{\baselinestretch}{1}\normalsize
L. 2/4: 't' is set up as referencing any fundamental line such as must 
exist when (L. 3) two different areas are within a given area, which
areas 
have solely boundary points in common, and have a fundamental line in 
common; which line referenced by 't' (L. 4) has as end-points a non 
boundary point of each of the two different areas, is wholly within them 
and intersects their common boundary line, not in its end-points, in one 
and only one point (this last simply to aid metal clarity). L. 5/6: Any
two 
't' lines with a common end-point are stipulated to have no other point
in 
common.  L. 7: A number of 't' lines are then set out relative to five 
areas - merely alluded to by superscripts A,B,C,D,E, because they do not 
require formal identification, in order to aid comprehension. These 't' 
lines "connect" two base areas 'A' and 'B' to each of the other three.
L. 
8/9: By suitable additions a symmetrical triplet of lines, each line 
connecting the two base areas through one of the others, is formed, on
the 
pattern of Th. 2.14. By renaming, the situation is as given. L. 10: The
't' 
line "connecting" the two base areas directly is then assumed to be not 
within a$^1$. But then P$^E$ is within, and P$^C$ is not within, a ring 
consisting wholly of 'D', 'A' and 'B' points and hence by Th. 2.13 and
2.15 
a 't' line "connecting" 'E' and 'C' cannot exist (L. 11). Similarly if
the 
'w' line is within a$^1$, the same problem arises between 'D', and 'E'
or 
'C' (L. 12/13). By symmetry all situations have now been examined so 't' 
lines cannot exist each "connecting" one of every pair out of five
areas, 
which would have to be the case if they each had with every other, a 
boundary line section in common.

\end{sffamily}
\end{document}

\begin{sffamily}
\large
NOTES

\normalsize
The abbreviations used to refer to different subsections in the paper
are as 
follows: D. = Definition; Cn. = 
Convention; In. = Interpretation; Th. = Theorem; Cr. = Corollary. In the 
third columns of proofs It. stands for iteration; Rn. for rename; Sim.
for 
similarly or similarity; Sym. for symmetry. In the commentary following 
proofs L. = Line. 

Cn. 1.1: The relationship between the objects referenced 
by a variable and that variable may be thought of as roughly similar to
the 
relationship between a mathematical variable and the values it takes. 
However no great weight should be placed on any such analogy. The usage
of 
the calculus demonstrates the nature of the concept being employed.  

The structure of proofs: Proofs are set out as numbered "lines" each
consisting 
of a statement which takes a truth value (such "lines" may occupy more
that 
one line on the page, but this is merely to facilitate the grouping of 
interlinked relationships, in order to aid comprehension.) The third
column 
of a proof indicates something of its structure, typically by listing 
previous lines as one or more premises on which the current line
depends, 
and often giving in addition the number of a relevant theorem which 
justifies that entailment, or other related information as alluded to in 
the abbreviations given above. Where a colon appears in this column it 
indicates one entailment embedded in another. However this third column
is 
not set out in accord with any set of rigorous rules from which the
exact 
structure of the proof can be extracted. The commentary which follows
most 
proofs should be referred to, in order to establish the proof's
structure. 
Typically omitted from the third column are the relationships indicated
in 
the commentary to Th. 1.1; and generally where a previous lines simply 
serves to define a variable or object, it is not listed.

Continuity: It may be asked what happens when you take an end-point,
say, away 
from a 
fundamental line ? What is left ? Certainly the rules of calculus allow 
that given \ul{AB}(x) one can assert the existence of $\ast$S(y) and 
y = [x - A]. 
One is then inclined to shift to the interpretive level and say, well 
what sort of "object" does 'y' represent then ? This however is an 
illegitimate move.  'y' cannot function as a variable referencing the 
points in a fundamental line because whatever end-point 'Z' one takes to 
accompany 'B' to form the end-points of that line, there always remains
a 
point referenced by 'y' between 'A' and 'Z' and hence not in that 
fundamental line. The set $\ast$S(y) then remains just that, a mere 
unstructured 
aggregation of points which cannot constitute a fundamental object of
the 
calculus. This is but one of the ways continuity, which is ruptured by
such 
a procedure, is expressed by the calculus.

Cn. 1.1: The syntax [A \map \ntt x] is not allowed because it can all
too 
easily lead to a drift into thinking of 
\ntt x as an \emph{object} \emph{of the calculus} where this is
illegitimate. 
For example consider two fundamental areas, with no boundary points in 
common, the first having within it the second. If 'a' references the
points 
of the "contained" area, to think of points "outside" the "contained"
area 
as mapping onto those points of the "containing" area which are "not a", 
tends to attract one into thinking of "not a" as referencing the points
of 
some kind of object or collection of objects of the calculus. But those 
points cannot constitute a well-formed array of the calculus analogously
to 
the case for a line noted above. Two fundamental areas suitably
conjoined 
and together having within them \emph{the boundary points} of the 
"contained" area, are necessary, in order to constitute the outer
annulus 
of the situation described. There thus seems to be two ways of
constructing 
the same situation, either as an annulus or as one area within another, 
which two different descriptions "compete for" the same set of points,
to 
form either the internal border of the annulus or the boundary of the 
contained area. This situation in fact appears to have a strong
resonance 
with the well known \emph{gestalt shift} between figure and ground: that
is 
the alternate perceptions of the same visual field as constituting
either 
an object lying on a surface, or an annular ring (against some
background). 
Whether any such parallel is of any significance however remains to be 
seen.

Cn. 2.1: If in order to avoid confusion it is necessary to indicate 
the seriate set relative to which a line is a seriating line one can do 
this for example by the syntax \ul{A-r-B}(x) where "r" is the variable 
referencing the constituting lines of the relevant seriate set. 

\vsp
\large
APPENDIX

\normalsize
The purpose of this appendix it to set out more rigorous proofs of Th.
1.1 and 
1.2 of the main text. In order to minimise repetition \emph{syntactic 
variables} will be employed. \emph{Syntactic variables} are variables
which 
are defined in a formal way in preliminary section, as referencing sets 
and/or objects all having the same, often complex, interlinking 
characteristics. In proofs these variables can then be directly alluded
to 
in order to assign those same sets of characteristics directly to
different 
objects without the need to list all those same characteristics every
time. 
These variables are indicated by the inclusion of '\%' as an element of 
their names, as in s\% and e$^C$\%(n). 

\renewcommand{\baselinestretch}{1.5}\small\normalsize
\begin{proof}
DEFINITION SV1:Syntactic variables s\%,e$^E$\%,e$^C$\%,p\%\\
1.\>[$\ast$S(n) \map s\%] \cd \{$\ast$S$_1$(m) \& [m \mm n]\}\\
2.\>\ul{J$^N$K}$^N$(j$^N$) \map n\\
3.\>\{[P \map j$^N$] \& \ntt[P = (J$^N$ \ori K$^N$)]\} \cd r\ntt\{[P 
\map j$^M$] \\
 \>\>\& [\ul{J$^M$K}$^M$(j$^M$) \map n] \& \ntt[j$^N$ \mm
j$^M$]\}\'\>1,2\\
4.\>P \map e$^E$\%:n\\
5.\>\{[P = (J$^M$ \ori K$^M$)] \& [\ul{J$^M$K}$^M$ (j$^M$) \map n]\\
 \>\>\& \ntt\{[P \map j$^L$] \& [\ul{J$^L$K}$^L$(j$^L$) \map n] \&
\ntt[j$^L$ 
 \mm j$^M$]\}\'\>1,4\\
6.\>P \map e$^C$\%:n\\
7.\>[P = (J$^M$ \ori K$^M$)] \& [P = (J$^L$ \ori K$^L$)]\\
 \>\>\&\{[\ul{J$^M$K}$^M$(j$^M$) \& \ul{J$^L$K}$^L$(j$^L$)]  \map n\} \& 
 \ntt[j$^L$ \mm j$^M$]\'\\
 \>\>\& \ntt\{[P \map j$^R$] \& [\ul{J$^R$K}$^R$(j$^R$) \map n]\'\\ 
 \>\>\& \ntt[j$^R$ \mm ( j$^M$ \ori j$^L$]\}\'\>1,6\\
8.\>\{[P = (J$^M$ \ori K$^M$)] \& [\ul{J$^M$K}$^M$(j$^M$) \map n]\}\\
 \>\>\cd [P \map (p$^E$\%:n \ori p$^C$\%:n)]\'\>1\\
9.\>[P \map p\%:n] \iif \{[\ul{J$^M$K}$^M$(j$^M$) \map n] \& [P \map
j$^M$]\}
\end{proof}

\renewcommand{\baselinestretch}{1}\small
Here the syntactic variable s\% is 
defined as referencing a set of fundamental lines, which have no points
in 
common except end-points, and whose each end-point is either, an element
of 
one and only one other of the lines of the set, or of no other such
line; 
and whose total set of subsumed points is referenced by p\%:n.

\renewcommand{\baselinestretch}{1.5}\normalsize
\begin{proof}
DEFINITION SV2: w\%\\
1.\>$\ast$S(n) \map w\%\\
2.\>$\ast$S(n) \map s\%\>\>1\\
3.\>P$^1$ \& P$^2$ \& \{[P$^1$ \& P$^1$] \map e$^E$\%:n\} \& 
\ntt[P$^1$ = P$^2$]l\\
 \>\>\& \ntt\{P$^3$ \& [P$^3$ \map e$^E$\%:n] \& \ntt\{P$^3$ = 
 (P$^1$ \ori P$^2$)]\}\'\>1
\end{proof}

\renewcommand{\baselinestretch}{1}\small
w\% references a set of fundamental lines which, in addition to being 
referenced by s\%, are such that  two and only two of the end-points of 
those lines are referenced by e$^E$\%. As a result such a set conforms
to 
the requirements of D. 0.2(b) clauses (ii)/(iv) for $\ast $S there.

\vsp
\normalsize
SV1/SV2 CONVERSE

If the entailments resulting from any set mapping onto either s\% 
and w\% are fulfilled by any set then that set maps onto those 
variables.

\vsp
THEOREM F1: Given two fundamental lines \ul{AB} and \ul{BC} which 
have no point in common except 'B', then the existence of a fundamental 
line \ul{AC}, constituted out of all and only the sum of the points of 
\ul{AB} and \ul{BC}, can be asserted.

\renewcommand{\baselinestretch}{1.5}\normalsize
\begin{proof}
1.\>\ul{AB}(x) \& \ul{BC}(y)\\
2.\>\ntt\{[P \map (x \& y)] \& \ntt[P = B]\}\\
3.\>R \& [R \map x] \& \ntt[R = (A \ori B]\\
4.\>$\ast$S$_($n) \& [$\ast$S$_($n) \map w\%] \& [R \map e$^C$\%:n]\\
 \>\>\& [(A \& B) \map e$^E$\%:n] \& [x \mm p\%:n]\'\>3,1,D.0.2,SV.1/2\\
5.\>z \mm (x + y)\\
6.\>$\ast$S$_($m) \& [$\ast$S(m) \map s\%] \& \{m \mm [n +
\ul{BC}(y)]\}\\
 \>\>\& [(B \& R) \map e$^C$\%:m] \& [(A \& C) \map e$^E$\%:m]\'\\
 \>\>\& [z \mm p\%:m]\'\>1,2,4,5)\\
7.\>$\ast$S(m) \map w\%6,SV.1/2,Conv.\\
8.\>$\ast$S$_($h) \& [$\ast$S(h) \map w\%] \& \{h \mm [\ul{AB}(x) + 
\ul{BC}(y)]\}\\
 \>\>\& [B  \map e$^C$\%:h] \& [(A \& C) \map e$^E$\%:h]\'\\
 \>\>\& [z \mm p\%:h]\'\>1,2,SV.1/2,Conv.\\
9.\>\ul{AC}(z) \& [z \mm (x + y)]\>\>7,It(R),Sm,8,D.0.2
\end{proof}

\renewcommand{\baselinestretch}{1}\small
L. 9: The result of L. 7 occurs for every 'R' in 'x' and similarly for
every 
non end-point of 'y'. L. 8 gives 
the same result for 'B'. Thus the requirement of a seriate set for every 
non end-point of 'z' is fulfilled. Here clause (v) of D. 0.2(b) is
ignored. 
In fact this clause is not strictly necessary to the definition of a 
seriate set and is only included in D. 0.2 to allow the intuitive syntax
in 
the introductory proofs of Th. 1.1 \& 1.2 of the main text, to be
directly 
employed. Clearly Th. F1 still goes through when that clause is ignored
- 
any such ring as it precludes simply existing before and after. Then, by 
application of Th. F1 itself and the axiom of stability of D. 0.4 the 
requirement of that clause can be proved. Intuitively this may seem
wrong, 
as if the truth of clause (v) must be proved before Th. F1 can be
proven. 
However here as elsewhere it is the entailments of the calculus which
must 
be attended to, not the interpretive situation where we may feel
compelled 
to think that the intuitive syntax used in the main text must be
justified 
before any such addition can occur.

\vsp
\normalsize
THEOREM F2: Given a non end-point C within a fundamental line \ul{AB}
then the 
existence of two fundamental 
lines \ul{AC} and \ul{CB} within \ul{AB} can be asserted, having C as 
their sole common point, and being constituted by all and only the same
set 
of points as constitute \ul{AB}.

\renewcommand{\baselinestretch}{1.5}\small\normalsize
\begin{proof}
1.\>\ul{AB}(x)\\
2.\>[C \map x] \& \ntt[C = (A \ori B)]\\
3.\>$\ast$S$_($n) \& [$\ast$S(n) \map w\%] \& [(A \& B) \map
e$^E$\%:n]\\
 \>\spc \& \{[P \map t] \iif [P \map e$^C$\%:n]\} \& [C \map $^t$]\\
 \>\spc \& [p\%:n \mm x]\>\>1,2,D.0.2,SV.1/2\\
4.\>\{[$\ast$S(j) \& $\ast$S(k)] \map w\%\}\\
 \>\spc \& [(j + k) \mm n] \& \ntt\{S$_1$! \& [S$_1$! \map (j \& k)]\}\\
 \>\spc \& [(A \& C) \map e$^E$\%:j] \& [(C \& B) \map e$^E$\%:k]\\
 \>\spc \& \{[P \map t$^1$] \iif [P \map e$^C$\%:j]\}\\
 \>\spc \& \{[P \map t$^2$] \iif [P \map e$^C$\%:k]\}\\
 \>\spc \& \{[t - C] \mm [t$^1$ + t$^2$]\}\\
 \>\spc \& [(p\%:j + p\%:k) \mm x]\>\>3,SV.1/2,Conv.\\
5.\>[p\%:j \mm y] \& [p\%:k \map z]\\
6.\>\ul{AC}(y) \& \ul{CB}(z)\\
 \>\spc \& [x \mm (y + z)]\\
 \>\spc \& \ntt\{\ntt[P = C] \& [P \map (y \&
z)]\}\>\>4,Th.F1,It,5,D.0.4
\end{proof}

\renewcommand{\baselinestretch}{1}\small
L. 4 derives from L. 3 by the single change 
of uncoupling the two lines of 'n' of which 'C' is the end point. L. 6:
All 
the lines of 'j' are "conjoined" to form a single fundamental line by
Th. 
F1 iterated, and similarly with 'k'. If after any addition the lines of
'j' 
or 'k' a line from one had an end-point in common with a line from the 
other, which point was not 'C' then, applying Th. F1, D. 0.4 would be 
controverted, hence the separation of the two sets.
\newpage
\setlength{\unitlength}{1cm}
\begin{picture}(16,21)

\thicklines
\put(1.41,20.37){\line(1,0){9.18}}
\put(1.41,19.53){\line(1,0){9.18}}
\put(1.01,14.36){\line(1,0){3.42}}
\put(4.43,12.92){\line(1,0){1.83}}
\put(1.01,11.48){\line(1,0){5.25}}
\put(11.25,14.24){\line(1,0){3.68}}
\put(9.95,11.33){\line(1,0){2.1}}
\put(13.99,11.33){\line(1,0){1.58}}
\put(1.01,7.03){\line(1,0){4.86}}
\thinlines
\put(1.41,20.97){\line(1,0){9.18}}
\put(1.41,18.93){\line(1,0){9.18}}
\put(1.41,18.26){\line(1,0){10.1}}
\put(5.54,19.95){\line(1,0){.92}}
\put(1.01,13.64){\line(1,0){3.42}}
\put(1.01,12.2){\line(1,0){5.25}}
\put(5.71,13.64){\line(0,1){.22}}
\put(.23,12.92){\line(0,1){1.81}}
\put(.23,12.92){\line(1,0){.6}}
\put(1.2,12.92){\line(1,0){.41}}
\put(1.98,12.92){\line(1,0){.88}}
\put(.23,10.97){\line(0,1){1.23}}
\put(11.25,13.45){\line(1,0){3.68}}
\put(9.95,12.13){\line(1,0){2.1}}
\put(13.99,12.13){\line(1,0){1.58}}
\put(10.55,13.45){\line(0,1){1.09}}
\put(13.02,10.61){\line(0,1){1.53}}
\put(1.98,7.2){\line(0,1){.47}}
\put(4.61,7.2){\line(0,1){.47}}
\put(3.2,7.2){\line(0,1){.14}}
\put(1.01,7.2){\line(0,1){.88}}
\put(5.87,7.2){\line(0,1){.88}}
\put(1.41,20.97){\vector(0,-1){.43}}
\put(10.59,20.97){\vector(0,-1){.43}}
\put(1.41,18.26){\vector(0,1){1.09}}
\put(10.59,18.26){\vector(0,1){1.09}}
\put(2.33,18.93){\vector(0,1){.42}}
\put(3.7,18.93){\vector(0,1){.42}}
\put(6.46,18.93){\vector(0,1){.42}}
\put(7.84,18.93){\vector(0,1){.42}}
\put(1.41,19.7){\vector(0,1){.5}}
\put(1.41,20.2){\vector(0,-1){.5}}
\put(10.59,19.7){\vector(0,1){.5}}
\put(10.59,20.2){\vector(0,-1){.5}}
\put(11.51,20.37){\vector(-1,0){.74}}
\put(11.51,19.53){\vector(-1,0){.74}}
\put(6.46,19.7){\vector(0,1){.5}}
\put(6.46,20.2){\vector(0,-1){.5}}
\put(1.01,14.19){\vector(0,-1){.38}}
\put(4.43,14.19){\vector(0,-1){.38}}
\put(1.01,13.47){\vector(0,-1){1.11}}
\put(1.79,13.47){\vector(0,-1){1.11}}
\put(2.86,12.37){\vector(0,1){1.1}}
\put(2.86,13.47){\vector(0,-1){1.11}}
\put(4.43,13.47){\vector(0,-1){.39}}
\put(4.43,12.75){\vector(0,-1){.39}}
\put(6.26,12.75){\vector(0,-1){.39}}
\put(1.01,12.03){\vector(0,-1){.38}}
\put(6.26,12.03){\vector(0,-1){.38}}
\put(5.71,13.64){\vector(-1,0){1.1}}
\put(.23,12.2){\vector(1,0){.59}}
\put(11.25,14.07){\vector(0,-1){.46}}
\put(14.93,14.07){\vector(0,-1){.46}}
\put(9.95,11.96){\vector(0,-1){.47}}
\put(12.05,11.96){\vector(0,-1){.47}}
\put(13.99,11.96){\vector(0,-1){.47}}
\put(15.57,11.96){\vector(0,-1){.47}}
\put(13.35,14.07){\vector(0,-1){.46}}
\put(10.55,13.45){\vector(1,0){.45}}
\put(13.02,12.13){\vector(-1,0){.78}}
\put(4.61,7.51){\vector(-1,0){2.64}}
\put(1.98,7.51){\vector(1,0){2.63}}
\put(5.87,7.92){\vector(-1,0){4.87}}
\put(1.01,7.92){\vector(1,0){4.86}}
\put(1.41,20.37){\circle*{.12}}
\put(10.59,20.37){\circle*{.12}}
\put(1.41,19.53){\circle*{.12}}
\put(10.59,19.53){\circle*{.12}}
\put(2.33,19.53){\circle*{.12}}
\put(3.7,19.53){\circle*{.12}}
\put(6.46,19.53){\circle*{.12}}
\put(7.84,19.53){\circle*{.12}}
\put(6.46,20.37){\circle*{.12}}
\put(1.01,14.36){\circle*{.12}}
\put(4.43,14.36){\circle*{.12}}
\put(1.01,13.64){\circle*{.12}}
\put(4.43,13.64){\circle*{.12}}
\put(4.43,12.92){\circle*{.12}}
\put(6.26,12.92){\circle*{.12}}
\put(1.01,12.2){\circle*{.12}}
\put(6.26,12.2){\circle*{.12}}
\put(1.01,11.48){\circle*{.12}}
\put(6.26,11.48){\circle*{.12}}
\put(1.79,13.64){\circle*{.12}}
\put(2.86,13.64){\circle*{.12}}
\put(1.79,12.2){\circle*{.12}}
\put(2.86,12.2){\circle*{.12}}
\put(4.43,12.2){\circle*{.12}}
\put(11.25,14.24){\circle*{.12}}
\put(14.93,14.24){\circle*{.12}}
\put(11.25,13.45){\circle*{.12}}
\put(14.93,13.45){\circle*{.12}}
\put(9.95,12.13){\circle*{.12}}
\put(12.05,12.13){\circle*{.12}}
\put(13.99,12.13){\circle*{.12}}
\put(15.57,12.13){\circle*{.12}}
\put(9.95,11.33){\circle*{.12}}
\put(12.05,11.33){\circle*{.12}}
\put(13.99,11.33){\circle*{.12}}
\put(15.57,11.33){\circle*{.12}}
\put(13.35,14.24){\circle*{.12}}
\put(13.35,13.45){\circle*{.12}}
\put(11.67,13.45){\circle*{.12}}
\put(12.2,13.45){\circle*{.12}}
\put(12.72,13.45){\circle*{.12}}
\put(10.37,12.13){\circle*{.12}}
\put(10.9,12.13){\circle*{.12}}
\put(11.42,12.13){\circle*{.12}}
\put(14.01,13.45){\circle*{.12}}
\put(1.01,7.03){\circle*{.12}}
\put(5.87,7.03){\circle*{.12}}
\put(1.98,7.03){\circle*{.12}}
\put(3.2,7.03){\circle*{.12}}
\put(4.61,7.03){\circle*{.12}}
\LARGE
\put(3.25,16.55){The structure of a seriate set of points}
\put(3.38,9.41){Theorem 1.1}
\put(9.68,8.48){Theorem 1.2}
\put(.49,4.96){Theorem 1.6}
\normalsize
\put(11.75,20.29){Seriate set S!}
\put(11.75,19.45){Set of seriate sets *S}
\put(11.75,18.32){'E' objects of seriate}
\put(11.75,17.98){ sets of *S identical to}
\put(11.75,17.65){'E' objects of S!}
\put(3.98,21.12){'E' objects of seriate set S!}
\put(3.25,19.93){'I' object of S!}
\put(3.25,18.59){'E' objects of seriate sets of *S}
\put(5.34,14.67){Set of seriate}
\put(5.34,14.33){sets equivalent}
\put(5.34,13.99){to AB}
\put(-.12,10.53){Set of seriate}
\put(-.12,10.2){sets equivalent}
\put(-.12,9.86){to AC}
\put(-.12,14.96){'I' object of AB \& AC = T}
\put(.59,14.29){A}
\put(4.61,14.29){B}
\put(4.01,12.85){B}
\put(6.45,12.85){C}
\put(.59,11.41){A}
\put(6.45,11.41){C}
\put(2.48,14.5){x}
\put(5.11,13.06){y}
\put(3.4,11.62){z}
\put(9.75,15.15){Set of seriate sets}
\put(9.75,14.81){equivalent to AB}
\put(12.42,10.31){Set of seriate sets of}
\put(12.42,9.97){AB, added by Th. 1.1}
\put(12.42,9.63){to form AC}
\put(12.42,9.3){}
\put(10.83,14.17){A}
\put(15.11,14.17){B}
\put(9.53,11.27){A}
\put(12.24,11.27){C}
\put(13.57,11.27){C}
\put(15.75,11.27){A}
\put(12.41,14.38){x}
\put(10.77,11.49){y}
\put(14.55,11.47){z}
\put(13.48,14.39){C}
\put(6.05,6.9){B}
\put(1.57,6.73){H}
\put(3.27,6.73){J}
\put(4.79,6.73){K}
\put(3.32,8.02){x}
\put(3.3,7.57){t}
\put(2.59,7.15){r}
\put(3.9,7.15){s}
\put(2.59,6.43){r'}
\put(3.76,6.41){t'}
\put(.59,6.9){A}
\put(11.13,13.33){\vector(-1,-1){1.1}}
\put(11.55,13.33){\vector(-1,-1){1.1}}
\put(12.08,13.33){\vector(-1,-1){1.1}}
\put(12.6,13.33){\vector(-1,-1){1.1}}
\put(13.23,13.33){\vector(-1,-1){1.1}}
\put(13.43,13.31){\vector(1,-2){.52}}
\put(15.01,13.31){\vector(1,-2){.52}}
\thicklines
\put(2.59,6.88){\oval(1.22,.34)[b]}
\put(3.29,6.88){\oval(2.63,1.09)[b]}
\thicklines
\put(12.5,4.2){\oval(1.2,2.4)}
\thinlines
\put(12.5,4.8){\oval(1.9,1.9)[t]}
\put(12.5,4.8){\oval(2.6,2.6)[tl]}
\put(12.5,4.8){\oval(3.3,3.3)[tl]}
\put(12.5,4.8){\oval(4,4)[t]}
\put(12.5,3.6){\oval(1.9,1.9)[b]}
\put(12.5,3.6){\oval(2.6,2.6)[bl]}
\put(12.5,3.6){\oval(4,4)[b]}
\put(11.2,3.6){\line(0,1){1.2}}
\put(11.55,3.6){\line(0,1){1.2}}
\put(13.45,3.6){\line(0,1){1.2}}
\put(14.5,3.6){\line(0,1){1.2}}
\put(10.85,4.8){\vector(0,-1){1.2}}
\put(10.5,5.1){\vector(0,-1){.3}}
\put(10.85,5.1){\vector(0,-1){.3}}
\put(10.85,4.5){\vector(0,1){.3}}
\put(10.5,3.3){\vector(0,1){.3}}
\put(11.2,3.3){\vector(0,1){.3}}
\put(11.2,3.9){\vector(0,-1){.3}}
\put(12.2,6.45){\vector(1,0){.3}}
\put(12.2,6.1){\vector(1,0){.3}}
\put(12.2,5.75){\vector(1,0){.3}}
\put(12.8,5.75){\vector(-1,0){.3}}
\put(12.8,2.65){\vector(-1,0){.3}}
\put(12.2,2.65){\vector(1,0){.3}}
\put(12.2,2.3){\vector(1,0){.3}}
\put(12.5,5.4){\circle*{.12}}
\put(12.5,3){\circle*{.12}}
\put(11.9,4.8){\circle*{.12}}
\put(11.9,3.6){\circle*{.12}}
\put(10.38,4.8){\line(1,0){.59}}
\put(10.38,3.6){\line(1,0){.94}}
\put(12.5,5.63){\line(0,1){.94}}
\put(12.5,2.18){\line(0,1){.59}}
\put(12.5,5.05){A}
\put(12.5,3.15){B}
\put(12.05,4.65){G}
\put(12.05,3.5){F}
\put(10.65,4.2){f}
\put(10.7,5.75){e}
\put(11.25,4.8){c}
\put(11.6,4.2){a}
\put(11.05,2.65){d}
\put(12.9,4.2){z}
\put(13.25,4.2){b}
\put(14.25,4.2){g}
\thicklines
\put(3,2.2){\oval(2,1)}
\thinlines
\put(2.5,2.2){\oval(1.7,1.7)[tl]}
\put(2.5,2.2){\oval(2.4,2.4)[l]}
\put(2.5,2.2){\oval(3.1,3.1)[l]}
\put(3.5,2.2){\oval(1.7,1.7)[tr]}
\put(3.5,2.2){\oval(2.4,2.4)[r]}
\put(3.5,2.2){\oval(3.1,3.1)[r]}
\put(2.5,3.75){\line(1,0){1}}
\put(2.5,3.4){\line(1,0){1}}
\put(2.5,3.05){\line(1,0){1}}
\put(2.5,1){\line(1,0){1}}
\put(2.5,.65){\line(1,0){1}}
\put(1.3,2.5){\vector(0,-1){.3}}
\put(1.3,1.9){\vector(0,1){.3}}
\put(1.65,2.5){\vector(0,-1){.3}}
\put(4.35,2.5){\vector(0,-1){.3}}
\put(4.7,1.9){\vector(0,1){.3}}
\put(4.7,2.5){\vector(0,-1){.3}}
\put(2.7,3.05){\vector(1,0){.3}}
\put(3.3,3.05){\vector(-1,0){.3}}
\put(2.8,.65){\vector(-1,0){.3}}
\put(3.2,.65){\vector(1,0){.3}}
\put(3,2.7){\circle*{.12}}
\put(2,2.2){\circle*{.12}}
\put(4,2.2){\circle*{.12}}
\put(1.18,2.2){\line(1,0){.59}}
\put(4.23,2.2){\line(1,0){.59}}
\put(3,2.93){\line(0,1){.24}}
\put(2.1,2.1){A}
\put(3.7,2.1){B}
\put(3,2.4){C}
\put(3,3.45){a}
\put(2.5,3.1){c}
\put(3.5,3.1){d}
\put(3,1.05){b}
\put(3,.7){x}
\LARGE
\put(10,-.1){Theorem 1.8}
\put(3.5,-1){Theorem 1.9}
\end{picture}
\newpage
\newsavebox{\TWTR}\savebox{\TWTR}(0,0)[bl]{%
\multiput(2,2.5)(.32,0){10}{\line(0,1){2.79}}
\thicklines
\put(2,2.5){\line(0,1){2.79}}
\put(5.2,2.5){\line(0,1){2.79}}
\put(2.31,3.67){\line(6,1){2.24}}
\thinlines
\put(2,5.5){\line(0,1){1.3}}
\put(5.2,5.5){\line(0,1){1.3}}
\put(5.2,6.59){\vector(-1,0){3.2}}
\put(5,6.59){\vector(1,0){.2}}
\put(2.31,5.5){\line(0,1){.8}}
\put(4.56,5.5){\line(0,1){.8}}
\put(4.56,6.09){\vector(-1,0){2.24}}
\put(4.36,6.09){\vector(1,0){.2}}
\put(2.31,3.67){\circle*{.12}}
\put(2.31,2.5){\circle*{.12}}
\put(2.31,5.29){\circle*{.12}}
\put(4.56,4.04){\circle*{.12}}
\put(4.56,2.5){\circle*{.12}}
\put(4.56,5.29){\circle*{.12}}
\put(3.27,3.83){\circle*{.12}}
\put(2.31,1.49){\vector(0,1){.8}}
\put(2.31,1.5){\line(-1,0){.8}}
\put(3.27,1.49){\vector(0,1){.8}}
\put(3.27,1.5){\line(1,0){.4}}
\put(4.56,1.49){\vector(0,1){.8}}
\put(4.56,1.5){\line(1,0){.4}}
\put(3.5,6.64){x}
\put(3.39,6.14){x$^1$}
\put(2.36,5.39){A$^P$}
\put(4.61,5.39){A$^Q$}
\put(2.41,2.2){B$^P$}
\put(4.66,2.2){B$^Q$}
\put(2.12,3.67){P}
\put(4.56,3.72){Q}
\put(3.23,3.92){J}
\put(3.83,3.72){g}
\put(-3.54,-.31){\makebox(5,2)[rt]{x$^P$}}
\put(3.87,1.39){x$^N$}
\put(5.16,1.39){x$^Q$}
}

\newsavebox{\TWFR}\savebox{\TWFR}(0,0)[bl]{%
\multiput(2,2.501)(.32,0){10}{\line(0,1){2.79}}
\thicklines
\put(2,2.5){\line(0,1){2.79}}
\put(5.2,2.5){\line(0,1){2.79}}
\put(2,2.5){\line(1,0){3.2}}
\put(2,5.29){\line(1,0){3.2}}
\thinlines
\put(2,5.5){\line(0,1){1.3}}
\put(5.2,5.5){\line(0,1){1.3}}
\put(5.2,6.59){\vector(-1,0){3.2}}
\put(5,6.59){\vector(1,0){.2}}
\put(2,2.5){\circle*{.12}}
\put(2,5.29){\circle*{.12}}
\put(5.2,2.5){\circle*{.12}}
\put(5.2,5.29){\circle*{.12}}
\put(3.28,2.5){\circle*{.12}}
\put(3.28,5.29){\circle*{.12}}
\put(3.5,6.64){x}
\put(1.65,3.8){x$^1$}
\put(5.29,3.8){x$^2$}
\put(1.5,5.17){A$^1$}
\put(1.5,2.37){B$^1$}
\put(5.33,5.17){A$^2$}
\put(5.33,2.37){B$^2$}
\put(3.08,5.45){A$^N$}
\put(3.28,2.12){B$^N$}
\put(2.88,3.5){\vector(1,0){.4}}
\put(2.88,3.5){\line(0,-1){1.7}}
\put(2.83,1.55){x$^N$}
\put(4.21,5.89){\vector(-1,-3){.2}}
\put(4.06,5.87){f}
\put(4.32,1.89){\vector(-1,3){.2}}
\put(4.41,1.61){g}
}
\newsavebox{\TWFV}\savebox{\TWFV}(0,0)[bl]{%
\multiput(2,2.5)(.32,0){10}{\line(0,1){2.79}}
\thicklines
\put(2,2.5){\line(0,1){2.79}}
\put(5.2,2.5){\line(0,1){2.79}}
\put(2.67,2.5){\line(1,0){2.11}}
\put(2.47,5.29){\line(1,0){1.43}}
\thinlines
\put(2,2.5){\line(1,0){3.2}}
\put(2,5.29){\line(1,0){3.2}}
\put(2,5.64){\line(0,1){.4}}
\put(5.2,5.64){\line(0,1){.4}}
\put(5.2,5.84){\vector(-1,0){3.2}}
\put(5,5.84){\vector(1,0){.2}}
\put(2,2.14){\line(0,-1){.4}}
\put(5.2,2.14){\line(0,-1){.4}}
\put(5.2,1.94){\vector(-1,0){3.2}}
\put(5,1.94){\vector(1,0){.2}}
\put(1.64,5.29){\line(-1,0){.4}}
\put(1.64,2.5){\line(-1,0){.4}}
\put(1.44,5.29){\vector(0,-1){2.79}}
\put(1.44,5.09){\vector(0,1){.2}}
\put(5.54,5.29){\line(1,0){.4}}
\put(5.54,2.5){\line(1,0){.4}}
\put(5.75,5.29){\vector(0,-1){2.79}}
\put(5.75,5.09){\vector(0,1){.2}}
\put(2.47,6.09){\line(0,1){.4}}
\put(3.91,6.09){\line(0,1){.4}}
\put(3.91,6.29){\vector(-1,0){1.43}}
\put(3.71,6.29){\vector(1,0){.2}}
\put(2.67,1.69){\line(0,-1){.4}}
\put(4.78,1.69){\line(0,-1){.4}}
\put(4.78,1.5){\vector(-1,0){2.11}}
\put(4.58,1.5){\vector(1,0){.2}}
\put(1,2.29){\line(0,1){3.19}}
\put(6.2,2.29){\line(0,1){3.19}}
\put(1.8,6.29){\vector(1,0){.67}}
\put(5.4,6.29){\vector(-1,0){1.48}}
\put(1.8,1.5){\vector(1,0){.87}}
\put(5.4,1.5){\vector(-1,0){.61}}
\put(1,5.5){\line(1,1){.8}}
\put(1,2.29){\line(1,-1){.8}}
\put(6.2,5.5){\line(-1,1){.8}}
\put(6.2,2.29){\line(-1,-1){.8}}
\put(4.62,5.29){\line(-1,-1){.8}}
\put(3.82,4.5){\line(-1,0){1.82}}
\put(2,3.62){\line(1,0){3.2}}
\put(4.62,5.29){\circle*{.12}}
\put(2,4.5){\circle*{.12}}
\put(2,3.62){\circle*{.12}}
\put(5.2,3.62){\circle*{.12}}
\put(3,4.5){\circle*{.12}}
\put(4,3.62){\circle*{.12}}
\put(5.2,3){\circle*{.12}}
\put(4.62,5.39){K$^4$}
\put(5.25,3.62){K$^3$}
\put(5.25,2.89){T$^3$}
\put(-3.05,2.5){\makebox(5,2)[rt]{J$^4$}}
\put(-3.05,1.62){\makebox(5,2)[rt]{J$^3$}}
\put(-1.8,2.45){\makebox(5,2)[rt]{T$^2$}}
\put(-.8,1.57){\makebox(5,2)[rt]{T$^1$}}
\put(3.62,4.5){y$^4$}
\put(2.95,3.62){y$^3$}
\put(4.62,4.12){a}
\put(2,2.5){\circle*{.12}}
\put(2,5.29){\circle*{.12}}
\put(5.2,2.5){\circle*{.12}}
\put(5.2,5.29){\circle*{.12}}
\put(2.47,5.29){\circle*{.12}}
\put(3.91,5.29){\circle*{.12}}
\put(2.67,2.5){\circle*{.12}}
\put(4.78,2.5){\circle*{.12}}
\put(-3.05,3.64){\makebox(5,2)[rt]{A$^1$}}
\put(-3.05,.5){\makebox(5,2)[rt]{B$^1$}}
\put(5.33,5.34){A$^2$}
\put(5.29,2.2){B$^2$}
\put(-4.05,2.17){\makebox(5,2)[rt]{h}}
\put(6.29,3.62){w}
\put(-3.58,2.04){\makebox(5,2)[rt]{x$^1$}}
\put(5.79,4.32){x$^2$}
\put(3.2,6.45){y$^1$}
\put(3.6,5.91){f}
\put(-1.28,-.56){\makebox(5,2)[rt]{y$^2$}}
\put(-1.4,-.11){\makebox(5,2)[rt]{g}}
\put(2.47,5.39){J$^1$}
\put(3.71,5.39){K$^1$}
\put(-2.33,.43){\makebox(5,2)[rt]{J$^2$}}
\put(-.22,.43){\makebox(5,2)[rt]{K$^2$}}
}
\newsavebox{\TWSVN}\savebox{\TWSVN}(0,0)[bl]{%
\thicklines
\put(4.87,4.87){\line(-3,-1){1.5}}
\put(3.37,4.37){\line(0,-1){.2}}
\put(3.37,4.17){\line(1,0){1}}
\put(4.37,4.17){\line(0,-1){1.19}}
\put(4.37,2.97){\line(-3,1){1.25}}
\put(3.12,3.39){\line(3,1){1}}
\put(4.12,3.72){\line(0,1){.25}}
\put(4.12,3.97){\line(-1,0){1.8}}
\put(2.33,3.97){\line(0,-1){.59}}
\put(2.33,3.37){\line(1,-1){.4}}
\put(2.72,2.97){\line(1,0){.4}}
\thinlines
\multiput(2,2.5)(.32,0){10}{\line(0,1){2.79}}
\thicklines
\put(2,2.5){\line(0,1){2.79}}
\put(5.2,2.5){\line(0,1){2.79}}
\put(2.54,2.5){\line(0,1){2.79}}
\put(3.53,2.5){\line(0,1){2.79}}
\put(4.49,2.5){\line(0,1){2.79}}
\thinlines
\put(2,5.39){\line(0,1){.9}}
\put(5.2,5.39){\line(0,1){.9}}
\put(5.2,6.09){\vector(-1,0){3.2}}
\put(5,6.09){\vector(1,0){.2}}
\put(3.5,5.59){\vector(-1,0){1.5}}
\put(3.3,5.59){\vector(1,0){.2}}
\put(5.16,5.59){\vector(-1,0){1.6}}
\put(4.96,5.59){\vector(1,0){.2}}
\put(2.54,1.79){\vector(0,1){.59}}
\put(2.54,1.8){\line(-1,0){.59}}
\put(-3.01,0){\makebox(5,2)[rt]{x$^P$}}
\put(3.53,1.79){\vector(0,1){.59}}
\put(3.53,1.8){\line(-1,0){.34}}
\put(-1.77,0){\makebox(5,2)[rt]{x$^N$}}
\put(4.49,1.79){\vector(0,1){.59}}
\put(4.49,1.8){\line(1,0){.4}}
\put(5.04,1.69){x$^Q$}
\put(3.12,2.97){\circle*{.12}}
\put(4.87,4.87){\circle*{.12}}
\put(2.54,3.97){\circle*{.12}}
\put(4.49,4.75){\circle*{.12}}
\put(4.87,4.92){Q}
\put(4.54,4.45){Q$^1$}
\put(2.64,4.07){P$^1$}
\put(3.12,2.67){P}
\put(3.5,6.14){x}
\put(2.76,5.64){f}
\put(4.36,5.7){g}
\put(-3.05,2.17){\makebox(5,2)[rt]{x$^1$}}
\put(5.29,3.89){x$^2$}
}
\newsavebox{\TWTN}\savebox{\TWTN}(0,0)[bl]{%
\thicklines
\put(4.24,7.29){\line(0,-1){.55}}
\put(4.24,6.75){\line(-1,0){.9}}
\put(3.34,6.75){\line(-3,2){.5}}
\put(2.84,7.08){\line(-2,-3){.59}}
\put(2.24,6.18){\line(-1,0){.5}}
\put(1.74,6.18){\line(0,-1){2.1}}
\put(1.74,4.08){\line(1,0){.8}}
\put(2.54,4.08){\line(1,-3){.19}}
\put(2.72,3.51){\line(1,0){1.6}}
\put(4.33,3.51){\line(-2,3){.26}}
\put(4.07,3.9){\line(1,6){.29}}
\put(4.37,5.7){\line(-1,0){.5}}
\put(3.87,5.7){\line(0,-1){.2}}
\put(3.87,5.5){\line(-1,0){.94}}
\put(2.91,5.5){\line(0,-1){.17}}
\put(2.6,4.55){\line(-1,0){.14}}
\put(2.45,4.55){\line(0,1){1.3}}
\put(2.45,5.85){\line(3,1){2.2}}
\put(4.65,6.59){\line(0,-1){3.39}}
\put(1.24,7.29){\circle{.11}}
\put(4.95,7.29){\circle{.11}}
\put(4.24,7.29){\circle*{.12}}
\put(1.24,6.75){\circle{.11}}
\put(4.95,6.75){\circle{.11}}
\put(4.24,6.75){\circle*{.12}}
\put(2.61,6.75){\circle{.11}}
\put(1.24,6.75){\circle{.11}}
\put(4.95,6.75){\circle{.11}}
\put(3.34,6.75){\circle*{.12}}
\put(2.61,6.75){\circle{.11}}
\put(1.24,7.08){\circle{.11}}
\put(4.95,7.08){\circle{.11}}
\put(2.84,7.08){\circle*{.12}}
\put(4.24,7.08){\circle{.11}}
\put(1.24,6.18){\circle{.11}}
\put(4.95,6.18){\circle{.11}}
\put(2.24,6.18){\circle*{.12}}
\put(3.43,6.18){\circle{.11}}
\put(4.65,6.18){\circle{.11}}
\put(1.24,6.18){\circle{.11}}
\put(4.95,6.18){\circle{.11}}
\put(1.74,6.18){\circle*{.12}}
\put(3.43,6.18){\circle{.11}}
\put(4.65,6.18){\circle{.11}}
\put(1.24,4.08){\circle{.11}}
\put(4.95,4.08){\circle{.11}}
\put(1.74,4.08){\circle*{.12}}
\put(4.1,4.08){\circle{.11}}
\put(4.65,4.08){\circle{.11}}
\put(1.24,4.08){\circle{.11}}
\put(4.95,4.08){\circle{.11}}
\put(2.54,4.08){\circle*{.12}}
\put(4.1,4.08){\circle{.11}}
\put(4.65,4.08){\circle{.11}}
\put(1.24,3.51){\circle{.11}}
\put(4.95,3.51){\circle{.11}}
\put(2.72,3.51){\circle*{.12}}
\put(4.65,3.51){\circle{.11}}
\put(1.24,3.51){\circle{.11}}
\put(4.95,3.51){\circle{.11}}
\put(4.33,3.51){\circle*{.12}}
\put(4.65,3.51){\circle{.11}}
\put(1.24,3.9){\circle{.11}}
\put(4.95,3.9){\circle{.11}}
\put(4.07,3.9){\circle*{.12}}
\put(2.6,3.9){\circle{.11}}
\put(4.65,3.9){\circle{.11}}
\put(1.24,5.7){\circle{.11}}
\put(4.95,5.7){\circle{.11}}
\put(4.37,5.7){\circle*{.12}}
\put(1.74,5.7){\circle{.11}}
\put(2.45,5.7){\circle{.11}}
\put(4.65,5.7){\circle{.11}}
\put(1.24,5.7){\circle{.11}}
\put(4.95,5.7){\circle{.11}}
\put(3.87,5.7){\circle*{.12}}
\put(1.74,5.7){\circle{.11}}
\put(2.45,5.7){\circle{.11}}
\put(4.65,5.7){\circle{.11}}
\put(1.24,5.5){\circle{.11}}
\put(4.95,5.5){\circle{.11}}
\put(3.87,5.5){\circle*{.12}}
\put(1.74,5.5){\circle{.11}}
\put(4.33,5.5){\circle{.11}}
\put(2.45,5.5){\circle{.11}}
\put(4.65,5.5){\circle{.11}}
\put(1.24,5.5){\circle{.11}}
\put(4.95,5.5){\circle{.11}}
\put(2.91,5.5){\circle*{.12}}
\put(1.74,5.5){\circle{.11}}
\put(4.33,5.5){\circle{.11}}
\put(2.45,5.5){\circle{.11}}
\put(4.65,5.5){\circle{.11}}
\put(1.24,4.55){\circle{.11}}
\put(4.95,4.55){\circle{.11}}
\put(1.74,4.55){\circle{.11}}
\put(4.17,4.55){\circle{.11}}
\put(4.65,4.55){\circle{.11}}
\put(1.24,4.55){\circle{.11}}
\put(4.95,4.55){\circle{.11}}
\put(2.45,4.55){\circle*{.12}}
\put(1.74,4.55){\circle{.11}}
\put(4.17,4.55){\circle{.11}}
\put(4.65,4.55){\circle{.11}}
\put(1.24,5.85){\circle{.11}}
\put(4.95,5.85){\circle{.11}}
\put(2.45,5.85){\circle*{.12}}
\put(1.74,5.85){\circle{.11}}
\put(4.65,5.85){\circle{.11}}
\put(1.24,6.59){\circle{.11}}
\put(4.95,6.59){\circle{.11}}
\put(4.65,6.59){\circle*{.12}}
\put(2.51,6.59){\circle{.11}}
\put(1.24,3.19){\circle{.11}}
\put(4.95,3.19){\circle{.11}}
\put(4.65,3.19){\circle*{.12}}
\put(1.24,4.75){\circle{.11}}
\put(4.95,4.75){\circle{.11}}
\put(1.74,4.75){\circle{.11}}
\put(4.21,4.75){\circle{.11}}
\put(2.45,4.75){\circle{.11}}
\put(4.65,4.75){\circle{.11}}
\put(1.24,4.75){\circle{.11}}
\put(4.95,4.75){\circle{.11}}
\put(1.74,4.75){\circle{.11}}
\put(4.21,4.75){\circle{.11}}
\put(2.45,4.75){\circle{.11}}
\put(4.65,4.75){\circle{.11}}
\put(1.24,4.95){\circle{.11}}
\put(4.95,4.95){\circle{.11}}
\put(1.74,4.95){\circle{.11}}
\put(4.24,4.95){\circle{.11}}
\put(2.45,4.95){\circle{.11}}
\put(4.65,4.95){\circle{.11}}
\thinlines
\put(3.21,4.75){\circle{1.3}}
\put(1.24,3.19){\line(0,1){4.1}}
\put(4.95,3.19){\line(0,1){4.1}}
\put(1.24,7.29){\line(1,0){3.71}}
\put(1.24,6.75){\line(1,0){3.71}}
\put(1.24,7.08){\line(1,0){3.71}}
\put(1.24,6.18){\line(1,0){3.71}}
\put(1.24,4.08){\line(1,0){3.71}}
\put(1.24,3.51){\line(1,0){3.71}}
\put(1.24,3.9){\line(1,0){3.71}}
\put(1.24,5.7){\line(1,0){3.71}}
\put(1.24,5.5){\line(1,0){3.71}}
\put(1.24,4.55){\line(1,0){1.36}}
\put(4.95,4.55){\line(-1,0){1.11}}
\put(1.24,5.85){\line(1,0){3.71}}
\put(1.24,6.59){\line(1,0){3.71}}
\put(1.24,3.19){\line(1,0){3.71}}
\put(1.24,4.75){\line(1,0){1.32}}
\put(4.95,4.75){\line(-1,0){1.14}}
\put(1.24,4.95){\line(1,0){1.36}}
\put(4.95,4.95){\line(-1,0){1.11}}
\put(4.14,7.39){J}
\put(4.55,2.82){K}
\put(5.75,6.99){\circle*{.12}}
\put(6.05,6.84){'f' points}
\put(5.75,6.54){\circle{.11}}
\put(6.05,6.39){'g' points also}
\put(6.05,6.09){ end-points of}
\put(6.05,5.79){'n' lines}
\put(3.21,4.1){\line(0,-1){1.5}}
}
\newsavebox{\sml}\savebox{\sml}(0,0)[bl]{%
\put(0,0){\oval(2.4738633753706,1.89)}
\thicklines
\put(-.6,.94){\line(0,-1){.94}}
\put(-.6,0){\line(1,0){.59}}
\put(0,0){\line(1,1){.4}}
\put(.4,.4){\line(1,-2){.4}}
\put(.8,-.41){\line(-1,0){2.03}}
\put(-.6,0){\circle*{.12}}
\put(0,0){\circle*{.12}}
\put(.4,.4){\circle*{.12}}
\put(.8,-.41){\circle*{.12}}
\put(-.6,.4){\circle{.11}}
\put(.59,0){\circle{.11}}
\thinlines
\put(-1.24,.4){\line(1,0){2.47}}
\put(-1.24,-.41){\line(1,0){2.47}}
\put(-1.24,0){\line(1,0){2.47}}
\put(-.5,.46){M$^1$}
\put(.46,.46){M$^2$}
\put(-1,-.3){N$^1$}
\put(.05,-.3){N$^2$}
\put(-1.81,.2){\vector(1,0){1.6}}
\put(1.92,.2){\vector(-1,0){1.6}}
\put(-6.97,-1.65){\makebox(5,2)[rt]{'s' area}}
\put(1.96,.1){'s' area}
}
\newsavebox{\TWAT}\savebox{\TWAT}(0,0)[bl]{%
\thicklines
\put(4.24,6.87){\line(3,-1){.72}}
\put(4.97,6.63){\line(-1,0){1}}
\put(3.97,6.63){\line(0,-1){.29}}
\put(3.97,6.33){\line(1,0){.29}}
\put(4.27,6.33){\line(2,-1){.5}}
\put(4.77,6.08){\line(-5,-6){.4}}
\put(4.37,5.6){\line(-1,0){.8}}
\put(3.57,5.6){\line(0,1){.9}}
\put(3.57,6.5){\line(-6,5){.5}}
\put(3.07,6.92){\line(-1,-2){.2}}
\put(2.87,6.52){\line(0,-1){.8}}
\put(2.87,5.72){\line(-5,-6){.3}}
\put(2.56,5.35){\line(3,-1){.59}}
\put(3.16,5.15){\line(0,-1){.29}}
\put(3.16,4.85){\line(-1,0){.8}}
\thicklines
\put(4.24,6.87){\circle{.13}}
\put(3.97,6.63){\circle{.13}}
\put(3.97,6.33){\circle{.13}}
\put(3.57,5.6){\circle{.13}}
\put(2.87,6.52){\circle{.13}}
\put(2.56,5.35){\circle{.13}}
\put(2.36,4.85){\circle{.13}}
\put(4.75,1.8){\circle{.13}}
\put(5,1.72){k points}
\put(4.97,6.63){\circle{.13}}
\put(4.77,6.08){\circle{.13}}
\put(3.57,6.5){\circle{.13}}
\put(2.87,5.72){\circle{.13}}
\put(3.16,5.15){\circle{.13}}
\put(3.16,4.85){\circle{.13}}
\put(2,4.5){\line(0,1){2.79}}
\put(5.52,4.5){\line(0,1){2.71}}
\put(4.4,3.9){\line(0,-1){1.3}}
\put(3.76,2.68){\line(-2,-1){1.32}}
\thinlines
\multiput(2,4.5)(.32,0){11}{\line(0,1){2.79}}
\put(2,4.09){\vector(1,0){2.35}}
\put(2.2,4.09){\vector(-1,0){.2}}
\put(5.52,4.09){\vector(-1,0){1.07}}
\put(5.32,4.09){\vector(1,0){.2}}
\put(3.1,4.14){s}
\put(4.9,4.14){d}
\put(2,4.29){\line(0,-1){.5}}
\put(2,3.4){\vector(0,1){.29}}
\put(1.75,3.1){S$_1$!$^1$}
\put(5.52,4.29){\line(0,-1){.5}}
\put(5.52,3.4){\vector(0,1){.29}}
\put(5.57,3.1){S$_1$!$^2$}
\put(4.29,6.91){Q}
\put(2.43,4.55){P}
\put(4.4,3.75){\circle*{.12}}
\put(4.5,3.65){F}
\put(4.4,3.25){\circle*{.12}}
\put(4.2,3){J}
\put(5.7,4.55){\vector(-1,-1){1.1}}
\put(6.1,5.85){\vector(-1,0){2.5}}
\put(6.1,4.95){\line(0,1){.9}}
\put(5.8,4.6){n$^B$ points} 
\put(4.4,2.75){\circle*{.12}}
\put(4.5,2.55){F}
\put(3.9,3.5){\circle*{.12}}
\put(3.6,3.5){F$^1$}
\put(4.9,3){\circle*{.12}}
\put(5,2.9){F$^1$}
\put(4.4,3.75){\line(-2,-1){.5}}
\put(3.9,3.5){\line(2,-1){1}}
\put(4.9,3){\line(-2,-1){.5}}

\put(3.1,3){\line(0,-1){1.3}}
\put(3.1,1.25){\vector(0,1){.4}}
\put(3.05,1.05){x$^J$}
\put(3.6,3){\line(0,-1){1.3}}
\put(3.6,1.25){\vector(0,1){.4}}
\put(3.55,1.05){x$^G$}
\put(2.6,2.3){\line(0,-1){.6}}
\put(2.6,1.25){\vector(0,1){.4}}
\put(2.5,1.05){x$^F$}

\put(3.1,2.85){\circle*{.12}}
\put(2.75,2.85){G$^1$}
\put(3.1,2.35){\circle*{.12}}
\put(3.15,2.15){J}
\put(1.48,2.89){\vector(3,-1){1.47}}
\put(1.48,2.85){\line(0,1){2.6}}
\put(-3.2,5.5){\makebox(5,2)[br]{s$^B$ points}}
\put(1.9,5.55){\vector(1,0){.82}}
\put(3.1,1.85){\circle*{.12}}
\put(3.2,1.7){F$^1$}
\put(3.6,2.6){\circle*{.12}}
\put(3.7,2.4){G} 
\put(2.6,2.1){\circle*{.12}}
\put(2.35,2.05){F}
\put(3.1,2.85){\line(2,-1){.5}}
\put(2.6,2.1){\line(2,-1){.5}}
}
\begin{picture}(16,21)
\put(0,15){\usebox{\TWTR}}
\put(9,15){\usebox{\TWFR}}
\put(0,6.9){\usebox{\TWFV}}
\put(9,7.1){\usebox{\TWSVN}}
\put(9.5,-1.5){\usebox{\TWAT}}
\put(0,-1.5){\usebox{\TWTN}}
\put(3,-1.5){\usebox{\sml}}
\LARGE
\put(2,14.9){Theorem 2.3}
\put(10.5,14.9){Theorem 2.4}
\put(2,7.5){Theorem 2.5}
\put(10.5,7.5){Theorem 2.8}
\put(7.5,-.2){Theorem 2.9}
\put(5.5,2.5){Theorem 2.11}

\end{picture}
\newpage
\newsavebox{\tpgn}\savebox{\tpgn}(0,0)[bl]{%
\thicklines
\put(1,2.35){\line(2,-3){.9}}
\put(1,2.35){\line(2,3){.9}}
\put(2.79,2.35){\line(-2,3){.9}}
\put(2.79,2.35){\line(-2,-3){.9}}
\put(1,2.35){\circle*{.12}}
\put(2.79,2.35){\circle*{.12}}
\put(1.89,1){\circle*{.12}}
\put(1.89,3.7){\circle*{.12}}
\put(.69,2.24){X}
\put(2.89,2.24){Y}
\put(1.84,3.8){P}
\put(1.89,.69){Q}
\put(1.3,3.04){x}
\put(2.39,1.44){y}
}\newsavebox{\tpa}\savebox{\tpa}(0,0)[bl]{%
\thicklines
\put(1.89,1){\line(0,1){2.7}}
\put(1.89,2.6){\circle*{.12}}
\put(1.59,2.5){Z}
\put(1.94,1.95){z}
}\newsavebox{\tpb}\savebox{\tpb}(0,0)[bl]{%
\thicklines
\put(1.89,3.7){\line(-1,0){1.5}}
\put(1.89,1){\line(-1,0){1.5}}
\put(.39,3.7){\line(0,-1){2.7}}
\put(.39,2.05){\circle*{.12}}
\put(.12,1.95){Z}
\put(.19,2.64){z}
}\newsavebox{\tpc}\savebox{\tpc}(0,0)[bl]{%
\thicklines
\put(1.89,3.7){\line(1,0){1.5}}
\put(1.89,1){\line(1,0){1.5}}
\put(3.39,3.7){\line(0,-1){2.7}}
\put(3.39,2.64){\circle*{.12}}
\put(3.5,2.5){Z}
\put(3.47,1.9){z}
}
\newsavebox{\ppt}\savebox{\ppt}(0,0)[bl]{%
\thicklines
\put(4,4.5){\line(0,1){3.29}}
\put(3,6){\line(1,0){2.29}}
\put(3,6){\circle*{.12}}
\put(2.7,5.88){F}
\put(5.29,6){\circle*{.12}}
\put(5.39,5.88){G}
\put(4,4.5){\circle*{.12}}
\put(3.79,7.89){P$^M$}
\put(4,7.79){\circle*{.12}}
\put(3.85,4.12){P$^N$}
\put(4,6){\circle*{.12}}
\put(3.75,5.67){T}
\put(3.85,6.4){j}
\put(4.2,6.04){k}
\put(4.45,4.79){a$^N$}
\put(4.4,7.39){a$^M$}
\put(2.89,3){Illustration of}
\put(2.89,2.5){\underline{P$^M$P$^N$} $\mapsto$ t}
\thinlines
\put(5.29,6){\line(1,1){.29}}
\put(4.79,6){\line(1,1){.79}}
\put(4.29,6){\line(1,1){1.29}}
\put(3.79,6){\line(1,1){1.79}}
\put(3.29,6){\line(1,1){2.09}}
\put(2.79,6){\line(1,1){2.09}}
\put(2.7,6.4){\line(1,1){1.69}}
\put(2.7,6.9){\line(1,1){1.19}}
\put(2.7,7.4){\line(1,1){.69}}
\put(5.29,6){\line(1,-1){.29}}
\put(4.79,6){\line(1,-1){.79}}
\put(4.29,6){\line(1,-1){1.29}}
\put(3.79,6){\line(1,-1){1.79}}
\put(3.29,6){\line(1,-1){1.79}}
\put(2.79,6){\line(1,-1){1.79}}
\put(2.7,5.59){\line(1,-1){1.39}}
\put(2.7,5.09){\line(1,-1){.89}}
\put(2.7,4.59){\line(1,-1){.39}}
}
\newsavebox{\fcp}\savebox{\fcp}(0,0)[bl]{%
\thicklines
\put(2.6,6){\line(1,-1){1.39}}
\put(2.6,6){\line(1,1){1.39}}
\put(5.4,6){\line(-1,1){1.39}}
\put(5.4,6){\line(-1,-1){1.39}}
\put(4,4.59){\circle*{.12}}
\put(4,7.4){\circle*{.12}}
\put(5.4,6){\circle*{.12}}
\put(5.5,5.88){P$^E$}
\put(3.95,7.5){P$^A$}
\put(3.89,4.15){P$^B$}
\put(4,6){\circle*{.12}}
\put(4.15,5.84){P$^D$}
}
\newsavebox{\lnn}\savebox{\lnn}(0,0)[bl]{%
\thicklines
\put(4,4.59){\line(0,1){2.79}}
\thinlines
\put(.75,6){\line(1,-1){3.25}}
\put(.75,6){\line(1,1){3.25}}
\put(7.25,6){\line(-1,1){3.25}}
\put(7.25,6){\line(-1,-1){3.25}}
\put(1.5,6){\vector(1,-1){1.94}}
\put(1.5,6){\vector(1,1){1.94}}
\put(6.5,6){\vector(-1,-1){1.94}}
\put(6.5,6){\vector(-1,1){1.94}}
\put(2.3,6){P$^C$}
\put(2.6,6){\circle*{.12}}
\put(4.09,6.34){j$^{AD}$}
\put(4.09,5.29){j$^{DB}$}
\put(3,6.79){j$^{AC}$}
\put(2.8,5){j$^{CB}$}
\put(4.79,6.79){j$^{AE}$}
\put(4.87,5.15){j$^{EB}$}
\put(5.9,5.2){w$^E$}
\put(2.35,7.15){w$^C$}
\put(1.35,7.39){\vector(0,-1){2.79}}
\put(1.35,7.19){\vector(0,1){.2}}
\put(1.14,6){w$^D$}
\put(5.5,4.25){\vector(-1,-1){.5}}
\put(5.7,4.45){\vector(1,1){.5}}
\put(5.71,4.17){h$^1$}
\put(1.25,3.5){Line 9}
}
\newsavebox{\lntn}\savebox{\lntn}(0,0)[bl]{%
\thicklines
\put(4,7.4){\line(1,0){2.23}}
\put(4,4.59){\line(1,0){2.23}}
\put(6.24,7.4){\line(0,-1){2.79}}
\thinlines
\put(4,4.59){\line(0,1){2.79}}
\put(1.6,6){\line(1,-1){2.1}}
\put(1.6,6){\line(1,1){2.1}}
\put(3.7,8.1){\line(1,0){3.23}}
\put(3.7,3.89){\line(1,0){3.23}}
\put(6.94,8.1){\line(0,-1){4.21}}
\put(6.3,5.29){w}
\put(2.3,6){P$^C$}
\put(2.6,6){\circle*{.12}}
\put(5.74,3.89){\vector(1,0){.5}}
\put(5.2,3.89){\vector(-1,0){.5}}
\put(5.2,3.94){h$^2$}
\put(2.6,3.1){Line 10}
}
\newsavebox{\lntw}\savebox{\lntw}(0,0)[bl]{%
\thicklines
\put(4,4.59){\line(0,1){2.79}}
\put(4,7.4){\line(-1,0){2.23}}
\put(4,4.59){\line(-1,0){2.23}}
\put(1.76,7.4){\line(0,-1){2.79}}
\thinlines
\put(3,6){\line(1,-1){1}}
\put(3,6){\line(1,1){1}}
\put(5,6){\line(-1,1){1}}
\put(5,6){\line(-1,-1){1}}
\put(2.6,5.5){w}
\put(1.86,5.9){P$^C$}
\put(1.76,6){\circle*{.12}}
\put(3.7,6.7){\vector(-1,-1){.5}}
\put(3.29,6.29){\vector(1,1){.5}}
\put(3.5,6.2){h$^3$}
\put(2.6,3.6){Line 12}
}
\begin{picture}(16,21)\thicklines
\put(-.5,17.5){\usebox{\tpgn}\usebox{\tpa}}
\put(6,17.5){\usebox{\tpgn}\usebox{\tpb}}
\put(11,17.5){\usebox{\tpgn}\usebox{\tpc}}
\put(-1.5,4.5){\usebox{\ppt}}
\put(5.5,3.5){\usebox{\fcp}\usebox{\lnn}}
\put(8.5,-2.5){\usebox{\fcp}\usebox{\lntn}}
\put(2.5,-2.5){\usebox{\fcp}\usebox{\lntw}}
\LARGE
\put(5,16){Theorem 2.14}
\put(3.5,-.5){A five country map}
\end{picture}
\end{sffamily}